\documentclass[12pt]{amsart}
\usepackage{amsthm,amsmath}
\usepackage{amssymb}

\newcommand{\al}{\alpha}
\newcommand{\be}{\beta}
\newcommand{\de}{\delta}
\newcommand{\ep}{\varepsilon}

\newcommand{\ee}{{\bf e}}
\newcommand\tc[2]{\theta\left[\begin{matrix}#1\\ #2\end{matrix}\right]}

\newcommand{\CC}{{\mathbb{C}}}

\newcommand{\EE}{{\mathbb{E}}}
\newcommand{\PP}{{\mathbb{P}}}
\newcommand{\QQ}{{\mathbb{Q}}}
\newcommand{\RR}{{\mathbb{R}}}
\newcommand{\ZZ}{{\mathbb{Z}}}

\newcommand{\calO}{{\mathcal O}}
\newcommand{\calH}{{\mathcal H}}
\newcommand{\calT}{{\mathcal T}}

\newcommand{\calA}{{\mathcal A}}

\newcommand{\calM}{{\mathcal M}}
\newcommand{\calJ}{{\mathcal J}}

\newcommand{\calG}{{\mathcal G}}
\newcommand{\calX}{{\mathcal X}}

\newcommand{\calR}{{\mathcal R}}
\newcommand{\calS}{{\mathcal S}}

\newcommand{\op}{\operatorname}
\newcommand{\Sat}{{\calA_g^{\op {Sat}}}}
\newcommand{\Vor}{{\calA_g^{\op {Vor}}}}
\newcommand{\Perf}{{\calA_g^{\op {Perf}}}}
\newcommand{\Perfl}{{\calA_g^{\op {Perf}}(\ell)}}

\newcommand{\Sp}{\op{Sp}}
\newcommand{\PSp}{\op{PSp}}

\newcommand{\Sym}{\op{Sym}}
\newcommand{\Sing}{\op{Sing}}

\newcommand{\tor}{\op{tor}}

\newcommand{\even}{\op{even}}

\def\diag{\operatorname{diag}}

\theoremstyle{plain}
\newtheorem{thm}{Theorem}[section]
\newtheorem{lm}[thm]{Lemma}
\newtheorem{prop}[thm]{Proposition}
\newtheorem{cor}[thm]{Corollary}

\newtheorem{conj}[thm]{Conjecture}

\theoremstyle{definition}

\newtheorem{rem}[thm]{Remark}

\begin{document}
\title[The locus of intermediate Jacobians]{The class of the locus of intermediate Jacobians of cubic threefolds}
\author{Samuel Grushevsky}
\address{Mathematics Department, Stony Brook University,
Stony Brook, NY 11790-3651, USA}
\email{sam@math.sunysb.edu}
\thanks{Research supported in part by National Science Foundation under the grant DMS-10-53313.}
\author{Klaus Hulek}
\address{Institut f\"ur Algebraische Geometrie, Leibniz Universit\"at Hannover, Welfengarten 1, 30060 Hannover, Germany}
\email{hulek@math.uni-hannover.de}
\thanks{Research is supported in part by DFG grants Hu-337/6-1 and Hu-337/6-2}

\begin{abstract}
We study the locus of intermediate Jacobians of cubic threefolds within the
moduli space $\calA_5$ of complex principally polarized abelian
fivefolds, and its generalization to arbitrary genus --- the locus of abelian varieties with a singular odd two-torsion point on the theta divisor. Assuming that this locus has expected codimension $g$ (which we show to be true for $g\le 5$, and conjecturally for any $g$),
we compute the class of this locus, and of is closure in the perfect cone toroidal compactification $\Perf$, in the Chow, homology, and the tautological ring.

We work out the cases of genus up to 5 in detail, obtaining explicit expressions for  the class of the
closures of $\calA_1\times\theta_{\rm null}$ in $\calA_4^{\op{Perf}}$, and for the class of the locus of intermediate Jacobians (together with the same locus of products) --- in $\calA_5^{\op{Perf}}$. Finally, we
obtain some results on the geometry of the boundary of the locus of intermediate Jacobians of cubic threefolds in $\calA_5^{\op{Perf}}$.

In the course of our computation we also deal
with various intersections of boundary divisors of a level toroidal compactification, which is of independent interest in understanding the cohomology and Chow rings of the moduli spaces.
\end{abstract}
\maketitle

\setcounter{section}{-1}
\section{Introduction}
\label{sec:intro}

The moduli spaces ${\mathcal M}_g$ of curves of genus $g$ and ${\mathcal A}_g$ of principally polarized abelian
varieties (ppav) of dimension $g$ are at the heart of algebraic geometry. Understanding their geometry includes the
question of computing the cohomology ring and the Chow ring of these varieties and their compactifications.

The investigation of the cohomology ring of ${\mathcal M}_g$ and its compactification $\overline{\calM_g}$,
the moduli space of stable curves, was pioneered by Mumford \cite{mumfordtowards}, further advanced by Faber's conjectures \cite{faberconjecture}, and is the contents of numerous
papers. The cohomology ring of $\calA_2$ was first computed by Igusa \cite{igugenus2}.
It was
Mumford \cite{mumfordtowards} who computed the cohomology and Chow ring of $\overline{\calM_2}$, or what is the same, of
the perfect cone compactification
$\calA_2^{\op{Perf}}$ (which coincides with both the Igusa and the second Voronoi compactification in this case).
In this paper he also sets up the general framework in which such computations are done nowadays,
in which the Grothendieck-Riemann-Roch theorem is a principal tool and the notion of tautological classes and
their relations play a crucial role. Mumford's conjecture on the stable cohomology of $\calM_g$ (as $g\to\infty$) was
famously proved recently by  Madsen and Weiss \cite{mawe}, and parts of Faber's conjecture on the tautological ring of $\calM_g$
were proved \cite{loo},\cite{ionel}, \cite{grva},\cite{faberonedim},\cite{boldsen}. However, many open questions remain, and in
particularly little is known about the suitably extended tautological ring of $\overline{\calM_g}$.

For the moduli space $\calA_g$ of ppav the situation is somewhat different: the tautological ring,
generated by the Chern classes $\lambda_i$
of the Hodge bundle, is known explicitly, and is quite easy to describe. Indeed, the only relation among the classes $\lambda_i$ on a suitable toroidal compactification $\overline{\calA_g}$ is
\begin{equation}\label{Rg}
(1+\lambda_1+\ldots+\lambda_g)(1-\lambda_1+\ldots+(-1)^g\lambda_g)=1,
\end{equation}
(while on $\calA_g$ one also has $\lambda_g=0$)
as proven by van der Geer \cite{vdgeercycles} in the Chow group with rational coefficients and by Esnault and Viehweg \cite{esvi} in the Chow ring. One can then compute the projections of various classes to the tautological ring --- for example, this was done by Faber \cite{faberalgorithms} for the locus of Jacobians, and by van der Geer \cite{vdgeercycles}, respectively Ekedahl and van der Geer \cite{ekvdgcycles} for the locus of products $[pt]\times\calA_{g-1}$. However, one sees that the classes of these loci do not lie in the tautological ring, and thus it is natural to define and study a larger subring in the Chow (or cohomology) ring that would contain such classes.
Defining such an extended tautological ring for the Deligne-Mumford compactification $\overline{\calM_g}$
or for a suitable toroidal compactification $\overline{\calA_g}$, or even the Satake compactification $\Sat$,
and understanding its structure, is a central
problem in understanding the intersection theory of compactified moduli spaces.

One can thus try to understand the structure of the Chow and homology of $\overline{\calM_g}$ and $\overline{\calA_g}$ in low genus.
The cohomology rings of $\calA_3$ and $\calA_3^{\op{Sat}}$ were computed by Hain \cite{hain}, the Chow rings of $\calA_3$  and $\calA_3^{\op{Perf}}$, with rational coefficients, were computed by
van der Geer \cite{vdgeercycles}, and finally the second-named author and Tommasi \cite{huto} computed the cohomology ring of $\calA_3^{\op{Perf}}$ and showed that it is equal to the Chow ring. Most of the cohomology of $\calA_4$ and its compactifications was computed by Tommasi and the second-named author \cite{huto}.

For higher genus, many questions remain open about the cohomology of $\overline{\calM_g}$ and $\overline{\calA_g}$. In particular, to the best of our knowledge for $g\ge 4$ {\em no} classes of geometrically meaningful cycles (of codimension higher than one --- the case of divisors is much easier) have been computed in the Chow ring of {\em any }toroidal
compactification of $\calA_g$, in particular the classes of $\calA_i \times \calA_{g-i}$ are unknown.

In this paper we concentrate on certain cycles of codimension $g$.
In genus $4$ this the class of the locus $\calA_1\times\calA_3$, but
our prime example is
the locus of intermediate Jacobians of cubic threefolds, its closure in $\calA_5^{\op{Perf}}$, and its generalization to arbitrary genus.
Assuming that such a generalized locus has codimension exactly $g$ in $\calA_g$ (for a detailed discussion of this condition see below), we compute the class of its closure in $\Perf$, and the projection of this class to the tautological ring.
In particular for $g=4$ and $g=5$ we determine the classes of the loci $\overline{\calA_1\times\theta_{\rm null}^{(3)}}\subset\calA_4^{\op{Perf}}$ and of the closure of the locus of intermediate Jacobians of cubic threefolds, together with $\overline{\calA_1\times\theta_{\rm null}^{(4)}}$, in $\calA_5^{\op{Perf}}$. We emphasize that we compute the classes of these loci precisely, not just their tautological parts.
Finally, we describe geometrically some strata of the boundary of the locus of intermediate Jacobians of cubic threefolds.

The locus of intermediate Jacobians of cubic threefolds and its closure is in manifold ways related to other interesting areas, foremost the theory of Prym varieties and the moduli space
$\calR_g$ of Prym covers, but also to the subvarieties of $\calM_g$
parameterizing curves which have a theta characteristic with a given number of sections,
and which have previously been studied by Harris \cite{hathetachar} and Teixidor i Bigas \cite{teixidor}. Moreover, our considerations have some bearing on the possible degenerations of Prym varieties, and are thus relevant for investigating geometrically meaningful
compactifications of the moduli space of cubic threefolds. We shall discuss these relations and some
connected questions in the subsequent section \ref{sec:statement}.

In the course of our computation we also study the combinatorics of intersections of boundary divisors of a level cover of $\Perf$, describing the classes of various geometric loci contained in the boundary of $\Perf$, which is of independent interest for understanding the structure of the Chow and cohomology rings.

Throughout the paper we work over the field of complex numbers.

\section*{Acknowledgements}
We would like to thank Gerard van der Geer for bringing the problem to our attention, and for insightful discussions on the geometry of the moduli space of abelian varieties. We are grateful to Valery Alexeev for conversations and explanations about the degenerations and polarizations on semi-abelic varieties and to Bill Fulton for
explanations about intersection theory.
The first author thanks Leibniz Universit\"at Hannover for hospitality during his visits there when some of this work was done. The second author is grateful to the Isaac Newton Institute in Cambridge, where this paper was completed. We are grateful to the referee for encouraging us to further explore the geometry of $\overline{IJ}$, which led to the creation of section 9.

\section{Statement of results and further outlook}
\label{sec:statement}
The intermediate Jacobians of cubic threefolds were studied in detail by Clemens and Griffiths \cite{clgr} who used them to prove that a non-singular cubic threefold is not rational. The locus
of intermediate Jacobians of cubic threefolds within the moduli space, which we denote $IJ^0\subset\calA_5$ is then a natural 10-dimensional geometrically defined quasi-projective subvariety of a 15-dimensional algebraic variety.
Similar to the Schottky problem for Jacobians of curves, it is interesting to try to describe the closure of this locus, which we denote $IJ\subset\calA_5$, by geometric or analytic conditions.

Clemens and Griffiths describe the geometry of the intermediate Jacobian in terms of the geometry of the threefold. We now review some of the beautiful geometric constructions associated to it. Indeed, consider the Prym map $p: \calR_6 \to \calA_5$, where $\calR_6$ is the space of
Prym curves of genus $6$ and recall that this map was compactified by Beauville \cite{beauville}.
The Prym map is generically finite of degree $27$ \cite{dosm}.
However, the situation is different over the locus
$IJ$ in $\calA_5$, where the general fiber has dimension $2$, see \cite{dofibers}.
Indeed, the fibers arise in the following way: if
$X$ is a cubic threefold, we consider a general line $l \subset X$. Projecting from this line gives $X$ birationally the
structure of a conic bundle over $\PP^2$. The discriminant curve $C \subset \PP^2$ of this
conic bundle is a smooth
quintic and hence of genus $6$. Every point on $C$ corresponds to a pair of different lines and
this gives rise to an
\'etale double cover $\tilde C \to C$, i.e. a point in $\calR_6$.

Using the Abel-Jacobi map, it was shown by Mumford \cite{mumfordprym},
see also \cite{beauvIJ} for detailed proofs,
that the Prym variety of this double cover is
isomorphic to the intermediate Jacobian $IJ(X)$. An analysis of the singularities of the
theta divisor, using the version of Riemann singularity theorem for Prym theta divisors, reveals that there is a unique triple point on the theta divisor. Since the
theta divisor is chosen to be symmetric, this triple point must then be an odd $2$-torsion point.
A geometric study of the Abel-Jacobi map further shows that the projectivized tangent cone to the theta divisor at this triple point is equal to the cubic threefold $X$, and in particular we see that $X$ can be recovered from its intermediate Jacobian. This is Mumford's \cite{mumfordprym} geometric proof of the
Torelli theorem for cubic threefolds, which is originally due to Clemens and Griffiths \cite{clgr}.

In view of the above, it is natural to ask whether all ppav with a triple point on the theta divisor are in fact intermediate
Jacobians of cubic threefolds. Using a detailed analysis of theta divisors of Prym varieties and degenerations,
Casalaina-Martin and Friedman \cite{cmfr}, \cite{casalaina2} answered this question in the positive: they showed that any indecomposable 5-dimensional ppav with a triple point on the theta divisor lies in the closure of the locus of intermediate Jacobians of cubic threefolds.
This statement can thus be interpreted as a geometric solution to the Schottky problem for intermediate Jacobians of cubic threefolds.

We denote by $A[2]$ the set of 2-torsion points on $A$, and add a subscript to distinguish the parity, and denote by $\calA_g(2)$ the full level two cover: the moduli of ppav together with a chosen symplectic basis for the group $A[2]$. We note that the intermediate Jacobian of a cubic threefold naturally comes with a chosen odd two-torsion point, but not naturally with a full level two structure. If we denote, following Donagi \cite{donagiscju}, $\calR\calA_g$ the moduli space of ppav together with one two-torsion point, which is a finite (but not Galois) cover of $\calA_g$, in turn covered by $\calA_g(2)$, we then get naturally the locus of intermediate Jacobians $\calR IJ\subset\calR\calA_5$. Note that this is the context in which Donagi \cite{donagiintermjac} showed that intermediate Jacobians of cubic threefolds are in the (big) Schottky-Jung locus. The question of understanding the geometry of $\calR\calA_5$ is of independent interest, and it seems very little is known.

\smallskip
The above characterization of $IJ\subset\calA_5$ can then be generalized to arbitrary genus to define the loci
$$
 I^{(g)}=\lbrace (A,\Theta)\in\calA_g\mid \Theta{\rm\ is\ singular\ at\ some\ point\
 m\in A[2]_{\rm odd}}\rbrace.
$$
The characterization of intermediate Jacobians then amounts to saying that within the locus of indecomposable ppav in $\calA_5$ we have $I^{(5)}=IJ$. However, the locus $I^{(5)}$ has an extra irreducible component consisting of decomposable ppav: in fact $I^{(5)}=IJ\cup (\calA_1\times\theta_{\rm null}^{(4)}$), where $\theta_{\rm null}^{(g)}\subset\calA_g$ denotes the theta-null divisor --- the locus of ppav for which the theta divisor is singular at an even 2-torsion point (recall that in our notation $IJ$ is the closure in $\calA_5$ of the locus $IJ^0$ of intermediate Jacobians).
The boundary of $IJ^0$, within $\calA_5$, and in the Satake compactification, was investigated by Casalaina-Martin and Laza \cite{cmla},
while an analytic description of the boundary of $I^{(5)}$ in the partial toroidal compactification was given by the first author and Salvati Manni in \cite{grsmconjectures}, using some of the ideas of their earlier work \cite{grsmodd1}.

More generally, it is natural to ask to describe the closures $\overline{IJ}$ and $\overline{I^{(g)}}$ of the loci $IJ$ and $I^{(g)}$ in some toroidal compactifications, and to study the possible degenerations of the cubic threefold itself, corresponding to various boundary strata of $\overline{IJ}$. A complete description of these could lead to a complete description of the geometry of a compactification of the moduli space of cubic threefolds. However, such a description would likely be extremely difficult, see \cite{alcato} for the complicated construction of this moduli space as a ball quotient. However, our results and computations in \cite{grhu2} give a way to describe geometrically some of the larger boundary strata --- see section \ref{degen} for more details.

The locus $I^{(g)}$ has expected codimension $g$ in $\calA_g$, and in \cite{grsmconjectures} it is conjectured that it is indeed of pure codimension $g$ (see section \ref{sec:notation} for a more detailed discussion).
It is thus natural to try to compute the class of $I^{(g)}$ in the Chow ring of $\calA_g$. The result in fact follows naturally from interpreting singularities of $\Theta$ at odd 2-torsion points as vanishing loci of gradients of theta functions, and thus as zero loci of certain vector-valued Siegel modular forms. The resulting expression is the content of our first main result:
\begin{thm}\label{thm:Igclass}
The virtual class of $I^{(g)}$ in $CH^g(\calA_g)$ is given by
$$
 [I^{(g)}]=2^{g-1}(2^g-1)\sum_{i=0}^{g}\lambda_{g-i}\left(\frac{\lambda_1}{2}\right)^i
$$
where $\lambda_i=c_i(\EE)$ are the Chern classes of the Hodge bundle on $\calA_g$.
\end{thm}
By virtual class we mean that if the codimension of the locus $I^{(g)}$ is $g$, as expected, then its class in the Chow group is given by the stated formula.
Note that this class lies in the tautological ring of the Chow ring generated by the classes $\lambda_i$. The locus $I^{(5)}$ is known to have expected codimension. This follows from the
work of Casalaina-Martin \cite{cmfr}, \cite{cmsurvey}, for a different proof see \cite{grhu2}.
Using the relations in $CH^*(\calA_g)$ it then follows that we have
\begin{cor}\label{cor:I5}
The class of $I^{(5)}$ in $CH^5(\calA_5)$ is equal to $$93\cdot(4\lambda_1^2\lambda_3+\lambda_1^5/2).$$
\end{cor}

\smallskip
Recall that the moduli space of ppav $\calA_g$ is not compact. In fact it admits many different toroidal compactifications. It is then natural to investigate the closure $\overline{I^{(g)}}$ in some toroidal compactification, and to compute its class there. The first steps in that direction were done in \cite{grsmconjectures} (see also \cite{cmla}) where the intersection of $\overline{I^{(g)}}$ with the boundary of the partial compactification were investigated.

The second (and much harder) main result of this paper is a computation of the class of $\overline{I^{(g)}}$ in the so-called perfect cone compactification $\Perf$.
Recall also that in genus $g\leq 3$ the perfect cone (also called first Voronoi) compactification
coincides with the second Voronoi, and the Igusa (also called central cone) compactification. The perfect cone
compactification has the advantage that its boundary divisor is irreducible. We also recall that
Shepherd-Barron \cite{shepherdbarron}
has shown that $\Perf$ is a canonical model in the sense of the minimal model program.
Computing the class $\overline{I^{(g)}}$ requires determining the vanishing behavior of the gradients of the theta functions on various loci of semiabelic varieties, and relies on the main results of our recent preprint \cite{grhu2}. We obtain the following
\begin{thm}\label{class}
For $g\le 5$ (and for any genus in which the statements of \cite[theorems 1.2 and 1.3]{grhu2} hold), we have the following expression for the class of $\overline{I^{(g)}}$ in $CH^g(\Perf)$:
\begin{equation}\label{Gclass}
[\overline{I^{(g)}}]=\frac{1}{N}\sum\limits_{m\in(\ZZ/2\ZZ)^{2g}_{\rm odd}}\sum\limits_{i=0}^g {p}_* \left( \lambda_{g-i}
\left( \frac{\lambda_1}{2} -\frac14 \sum\limits_{n\in Z_m} \de_n\right)^i \right)
\end{equation}
where $p :\Perf(2)\to\Perf$ is the level cover, $N=|\Sp(2g,\ZZ / 2\ZZ)|$ and
$Z_m$ is the set of pairs of non-zero vectors $\pm n \in (\ZZ/2\ZZ)^{2g}$ such that $m+n$ is even, and we recall that the irreducible components $\de_n$ of the boundary of $\Perf(2)$ correspond to non-zero elements of $(\ZZ/2\ZZ)^{2g}$.
\end{thm}
Here $(\ZZ/2\ZZ)^{2g}_{\rm odd}$ denotes the set of {\em odd} elements in $(\ZZ/2\ZZ)^{2g}$. Any such element
can be written in the form $m=(m_1,m_2)$ where $m_i \in (\ZZ/2\ZZ)^{g}$. We call $m$ odd if
the scalar product $m_1 \cdot m_2$ is $1$ and even otherwise.
At this point some words about intersection theory are in order. We denote by $CH^k(\Perf)$ the Chow group of cycles
of codimension $k$. We first note that this is always meant in the sense of the stack, in particular we are
free to
work with invariant classes on level covers and then take the pushforward to $\Perf$.
Secondly, we would like to point out that,
since the perfect cone decomposition contains non-basic cones for $g\geq 4$, the stack $\Perf$
is not smooth and hence there is no ring structure on $CH^*(\Perf)$.
However, we shall mostly be working with Chern classes of vector bundles which are elements in the
operational Chow ring, where we can multiply these classes and then take the cap product
with any cycle.

We shall also compute the projection of the class computed above to the {\em tautological ring}, i.e.~the ring generated by the Hodge classes:
recall that the tautological ring of a toroidal compactification $\calA_g^{\op{tor}}$ is the polynomial ring generated by the classes $\lambda_1,\ldots,\lambda_g$ subject to the one fundamental relation (\ref{Rg}).
This is
defined for every toroidal compactification of $\calA_g$ and independent of the chosen compactification, as is its
pushforward to the Satake compactification (see \cite{ekvdgcycles}).

The tautological ring is also defined for the open part $\calA_g$ where the extra relation $\lambda_g=0$ holds, which was shown by van der Geer \cite{vdgeercycles} in cohomology and by  Esnault and Viehweg \cite{esvi} in the Chow ring. We recall that the tautological ring
has a perfect pairing, and thus that there is a projection from the Chow ring to the tautological ring. For details we refer the reader to \cite{faberalgorithms}, \cite{vdgeercycles},\cite[section 3]{ekvdgcycles}.

\begin{thm}\label{thmtaut}
If \cite[theorems 1.2 and 1.3]{grhu2} hold in genus $g$ (in particular for any $g\le 5$)
the projection of the class
$[\overline{I^{(g)}}]$ to the tautological ring is equal to
$$
 [I^{(g)}]^{\op{taut}}=\frac{(-1)^{g-1}(g-1)!}{8\zeta(1-2g)}\lambda_g+
 2^{g-1}(2^g-1)\sum_{i=0}^{g}\lambda_{g-i}\left(\frac{\lambda_1}{2}\right)^i.
$$
\end{thm}
We will discuss the above statements, level covers, and their boundary components in detail in the following sections.

We apply the above theorems to compute the classes of $I^{(g)}$ and $\overline{I^{(g)}}$. for all $g\leq 5$. For $g=2$ these loci are empty, and we get zero as a valid consistency check for our computations. For genus 3 our results agree with the computation of the class  $[\overline{\Sym^3(\calA_1)}]\in CH^*({\calA_3}^{\op {Perf}})$
obtained by van der Geer in \cite{vdgeerchowa3}, and corrected in \cite{corr}. For genus 4 we compute the class of the locus $\calA_1\times\theta_{null}^{(3)}$ and of its closure in $\calA_4^{\op{Perf}}$. In particular, we confirm that the class of the open part of this locus lies in the tautological subring of $CH^*(\calA_4)$.
Finally, for genus 5 our results are completely new, and give a formula for the class of the locus of intermediate Jacobians of cubic threefolds together with the locus $\calA_1\times\theta_{\rm null}^{(4)}$, and also of the compactification. We
further compute the tautological projections and obtain for the closure  of the locus of
intermediate Jacobians of cubic threefolds $\overline{IJ}$:
\begin{prop}\label{IJtaut}
The projection of the class of the closure of the locus of intermediate Jacobians of cubic threefolds to the tautological ring is given by
$$
[\overline{IJ}]^{\op{taut}}=140\lambda_1^5 - 376 \lambda_1^2\lambda_3 + 848 \lambda_5.
$$
\end{prop}

\medskip
Except the calculations of the projections of various classes to the tautological ring (which depend on the existence of a perfect pairing and thus require working in the Chow ring with rational coefficients), all
our calculations hold in the Chow ring with integer coefficients.

At this point we would like to link our results to other areas and some open problems.
The first connection concerns the geometry of the Prym map,
in particular in genus $5$.

We have already noted that a general line on a cubic threefold $X$ gives rise to
a Prym variety which is isomorphic to the intermediate Jacobian of $X$. This shows that
the fiber of the Prym map $p: \calR_6 \to \calA_5$ has dimension at least $2$.
It was shown in \cite[4.6]{dofibers} that the fiber $p^{-1}(X)$ is indeed isomorphic to the Fano variety of $X$.

Thus the cycle $p^{-1}(IJ)\subset\calR_6$ has dimension $12$.
It is an interesting question to ask what the class of
this cycle, respectively its closure in Beauville's partial compactification $\calR'_g$ or the
compactification by stable curves $\overline{\calR}_g$, is.
As far as we know this problem is completely open.

One can also ask what happens to the other component of $I^{(5)}$,
namely $\calA_1 \times \theta_{\rm null}^{(4)}$
under pullback by the (compactified) Prym map $p': \calR_g' \to \calA_{g-1}$.
Indeed, the fiber of the compactified Prym map over a point of $\calA_1 \times \calA_4$
or $\calA_1 \times \theta_{\rm null}^{(4)}$ has not been considered in full detail in the literature.
It is easy to see that the dimension of $(p')^{-1}(\calA_1 \times \theta_{\rm null}^{(5)})$ is
at least (and presumably also equal to) $12$. Again, it would be interesting to know the class of this cycle and its closure.

The cycles $I^{(g)}$ also relate to interesting cycles on $\calM_g$. Recall the subvarieties
${\mathcal M}_g^k$ of $\calM_g$ of curves of genus $g$ having
a theta characteristic with at least $k+1$ sections and the same parity as $k+1$.
These were first introduced by Harris and studied by Harris \cite{hathetachar} and Teixidor i Bigas \cite{teixidor}.
In particular, it follows from their work that the codimension of the locus ${\mathcal M}_g^2$ in $\calM_g$ is $3$.
By the Riemann singularity theorem $I^{(g)} \cap \calM_g = {\mathcal M}_g^2$. Note that this intersection is
highly non-transversal; in particular its codimension is not $g$ and thus we cannot compute it by looking at
the top Chern class of the pullback of a rank $g$ vector bundle. It would be highly interesting to compute the classes of the cycles  ${\mathcal M}_g^k$ and
their closures.

\section{Gradients of theta functions}
\label{sec:notation}

We denote by $\calH_g$ the Siegel upper half space of genus $g$ --- the set of symmetric complex $g\times g$ matrices with positive-definite imaginary part. Recall that the symplectic group $\Sp(2g,\ZZ)$ acts on $\calH_g$ by $\gamma\circ\tau=(A\tau+B)(C\tau+D)^{-1}$.

We recall that the level subgroups of $\Gamma_g:=\Sp(2g,\ZZ)$ are defined as follows:
$$
 \Gamma_g(n):=\left\lbrace \gamma=\begin{pmatrix} A&B\\ C&D\end{pmatrix}\in\Gamma_g \right|\left. \gamma\equiv \begin{pmatrix} 1&0\\ 0&1\end{pmatrix}\mod n\right\rbrace
$$
$$
 \Gamma_g(n,2n):=\left\lbrace\gamma\in\Gamma_g(n)\mid \diag(A^tB)\equiv\diag(C^tD)\equiv 0
 \mod 2n\right\rbrace.
$$
The moduli space of ppav is then $\calA_g=\calH_g/\Gamma_g$, while the level moduli spaces $\calA_g(n):=\calH_g/\Gamma_g(n)$ and $\calA_g(n,2n):=\calH_g/\Gamma_g(n,2n)$ are finite covers of $\calA_g$.

We denote by $\theta(\tau,z)$ the Riemann theta function of $\tau\in\calH_g$ and $z\in\CC^g$
$$
 \theta(\tau,z):=\sum\limits_{n\in\ZZ^g}\ee (n^t\tau n/2 + n^tz)
$$
where for future use we denote $\ee(x):=\exp(2\pi i x)$ the exponential function.

For an abelian variety $A$, we denote by $A[2]$ the set of two-torsion points on it; as a group, $A[2]\cong (\ZZ/2\ZZ)^{2g}$.
Analytically the points in $A[2]$ can be labeled  $m=(\tau\ep+\de)/2$, where $\tau\in\calH_g$ projects to $A\in\calA_g$, and $\ep,\de\in (\ZZ/2\ZZ)^g$. For future use we denote $\sigma(m):=\ep\cdot\de\in\ZZ/2\ZZ$ and call it the parity of $m$.
Accordingly we
call $m$ even or odd depending on whether $\sigma(m)$ is 0 or 1, respectively.
This is equivalent to the point $m$ not lying (resp.~lying) on the theta divisor for a generic $\tau$ (i.e.~for a two-torsion point $m$
the function $\theta(\tau,m)$ is identically zero if and only if $m$ is odd).

For a point $m=(\tau\ep+\de)/2$ we denote the theta function with (half-integer) characteristic
$$
 \tc\ep\de(\tau,z):=\theta_m(\tau,z):=
$$
$$
=\sum\limits_{n\in\ZZ^g}\ee ((n+\ep/2)^t\tau (n+\ep/2)/2
 + (n+\ep/2)^t(z+\de/2)).
$$
As a function of $z$, the theta function with characteristic is even or odd depending on the parity of the characteristic. In particular, all theta constants (the values of theta functions with characteristics at $z=0$) with odd characteristics vanish identically.

Let $\pi:\calX_g\to\calA_g$ be the universal family, which exists over the stack, and let
$\EE:=\pi_*\Omega_{\calX_g/\calA_g}$ be the Hodge bundle.
The gradient
\begin{equation}\label{Fmdef}
  F_m:={\rm grad}_z\theta_m(\tau,z)|_{z=0}
\end{equation}
with respect to $z$ of the theta function vanishes identically in $\tau$ for even $m$, and is generically non-zero for $m$ odd.
This gradient is a vector-valued modular form for $\Gamma_g(4,8)$ for the
representation $\det^{\otimes 1/2}\otimes{\rm std}:\Gamma_g\to GL(g,\CC)$. In other words, we have
\begin{equation}\label{Fmbundle}
 F_m\in H^0(\calA_g(4,8),\det\EE^{\otimes 1/2}\otimes\EE)
\end{equation}
(see \cite{grsmodd1} for more details).

We recall from \cite[p.~50]{igusabook}, that
up to a simple exponential factor, the theta function with characteristic $m$ is equal to the
Riemann theta function shifted by the point $m$:
$$
 \theta(\tau,z+(\tau\ep+\de)/2)=\ee(-\ep^t\tau\ep/8-\ep^t\de/4-\ep^tz/2) \tc\ep\de(\tau,z).
$$
We can thus compute for an odd point $m\in A[2]$
\begin{equation}\label{fmdef}
 f_m(\tau):={\rm grad}_z\theta(\tau,z)|_{z=m}={\rm grad}_z\theta(\tau,z+(\tau\ep+\de)/2)_{z=0}
\end{equation}
$$
= \ee(-\ep^t\tau\ep/8-\ep^t\de/4-\ep^tz/2){\rm grad}_z\tc\ep\de(\tau,z)|_{z=0}
$$
$$
 =\ee(-\ep^t\tau\ep/8-\ep^t\de/4)F_m(\tau)
$$
since $\tc\ep\de(\tau,0)=0$ for the odd two-torsion point. Thus $f_m$ and $F_m$ differ by a
nowhere vanishing holomorphic function on $\calH_g$, and thus their zero loci are the same.
Moreover, the line bundle on $\calA_g(4,8)$ defined by the exponential factor is trivial
since it has a nowhere vanishing section. (In what follows, it will be crucial that this exponential factor, while non-vanishing on $\calA_g$, in fact vanishes on some irreducible components of the boundary of $\Perf(2)$.)
Thus we also have
$$
f_m\in H^0(\calA_g(4,8),\det\EE^{\otimes 1/2}\otimes\EE).
$$

The group $\Gamma_g(2)/\Gamma_g(4,8)$ acts on the gradients by certain signs, and thus the zero locus
\begin{equation}\label{Gmdef}
 G_m:=G_{\epsilon,\delta}
=\{\tau| F_m(\tau)=0\}
 =\{\tau| f_m(\tau)=0\}
\end{equation}
is a well-defined subvariety of $\calA_g(2)$.
Moreover, we note that to define $f_m$, we only need to choose one (odd) two-torsion point, and thus we have a well-defined zero locus $G_m$ of $f_m$ on $\calR\calA_g$. In principle one could then work with this non-Galois cover of $\calA_g$; however, the theory of toroidal compactifications is better developed for full level covers, and while not much is known, or could be now done for $\overline{\calR\calA_g}$, our results below are for $\overline{\calA_g(2)}$.

We refer to \cite{grsmodd1,grsmconjectures} for a
more detailed discussion of the properties of the
gradients of the theta function and further questions on loci of ppav with points of high multiplicity on the theta divisor.

Finally we denote by
$$
 I^{(g)}:=p(G_m)\subset\calA_g
$$
the locus of ppav for which {\it some} $F_m$ vanishes. Geometrically, as explained in the introduction, this is the locus of ppav having a singularity at an odd two-torsion point. We will omit the index $(g)$ when no confusion is possible. Note that it follows from the fact that $\Gamma_g$ permutes the $F_m$ that the
projection of $G_m$ to $\calA_g$ does not depend on $m$.

The multiplicity of the theta function for ppav of low dimension has been studied extensively. Recall that a ppav is called decomposable if it is a product of two lower-dimensional ppav.
For genus up to 4 it is known that no indecomposable ppav has a point of multiplicity 3 on the theta divisor \cite{cmsurvey}, and thus by
studying the multiplicity of points on the reducible theta divisors for decomposable ppav, we see that as schemes
$$
 I^{(3)}=\Sym^3(\calA_1)
$$
and
$$
 I^{(4)}=\calA_1\times\theta_{\rm null}^{(3)}
$$
where $\theta_{\rm null}^{(3)}$ denotes the theta-null divisor in $\calA_3$: the locus of those ppav for which there exists a point $m\in A[2]_{\rm even}$ lying on the theta divisor (or, equivalently, for which some theta constant vanishes).

It was recently shown in \cite{cmfr,casalaina2} that within the locus of indecomposable abelian 5-folds the locus $I^{(5)}$ coincides with the
locus $IJ$ of intermediate Jacobians of cubic threefolds, while in \cite{cmla} degenerations of intermediate Jacobians were studied. Combining this
with results of \cite{grsmconjectures} one gets scheme-theoretically
$$
 I^{(5)}=IJ\,\cup \,\calA_1\times\theta_{\rm null}^{(4)}.
$$
Recall again that $IJ$ denotes the closure in $\calA_5$ of the locus $IJ^0$ of intermediate Jacobians of cubic threefolds. In \cite{cmla} the boundary of $IJ$ was described, and in particular it was shown that within $\calA_5$ one has
$$
 IJ\setminus IJ^0=\overline{\calJ_5^h}\,\cup\,\calA_1\times(\overline{\calJ_4}\cap\theta_{\rm null}^{(4)}),
$$
where $\overline{\calJ_g}\subset\calA_g$ denotes the closure of the locus of Jacobians of curves, and $\overline{\calJ_g^h}\subset\calA_g$ denotes the closure of the locus of hyperelliptic Jacobians (see also \cite{grsmconjectures} for a discussion). In particular one sees that $\calA_1\times(\overline{\calJ_4}\cap\theta_{\rm null}^{(4)})$ is the intersection of the two irreducible components of the locus $I^{(5)}$.

\begin{rem}\label{remred}
The explicit descriptions of the loci $I^{(g)}$ given above for $g\le 5$ are a priori only set-theoretic. To make these descriptions scheme-theoretic, one would need
to ascertain that the gradients vanish with multiplicity one along each irreducible component of the locus. We do not know how to do accomplish this in general, and thus in principle theorem \ref{class} should be a priori interpreted as giving the class of the scheme $I^{(g)}$ which may have non-reduced components. However, for genus $g\le 5$ all the irreducible components of $I^{(g)}$ are known; in particular, it is known that they all intersect the boundary of the partial toroidal compactification, and our considerations in section \ref{degen} then show that each component intersects the boundary in a reduced scheme, and thus is itself reduced. We also note that
$\Perf$ is Cohen-Macaulay since toric varieties have this property. It then follows that $\overline{I^{(g)}}$ is also
Cohen-Macaulay provided it has codimension $g$, since a local complete intersection variety in a CM variety
is again CM. In particular $\overline{I^{(g)}}$ has no embedded components if it is of the correct codimension.
\end{rem}

In fact, in higher genus the locus $I^{(g)}$ is not well understood. Indeed, even the following question is open:
\begin{conj}\cite[Conjecture 1]{grsmconjectures}\label{codimconj}
The locus $I^{(g)}$ has pure codimension $g$ in $\calA_g$ for any $g$, and is reduced.
\end{conj}

Notice that since locally $I^{(g)}$ is given by the vanishing of $g$ partial derivatives
of the theta function at an odd two-torsion point,
we know that codimension of each its irreducible component is at most $g$. What the conjecture says
is thus that the vanishing of the partial derivatives imposes $g$ independent conditions.
In this paper we will concentrate on the case of $g\le 5$ (when this conjecture is known to hold (see \cite{cmsurvey} and \cite[Thm 1.2]{grhu2}). The above analytic description of $I^{(g)}$ allows us to prove our first main result.
\begin{proof}[Proof of theorem \ref{thm:Igclass}]
This is an immediate consequence of the fact that $I^{(g)}$ is the zero set of the sections $F_m$, resp. $f_m$
of the rank $g$ vector bundle $\EE \otimes \det \EE^{\otimes 1/2}$.
Provided $I^{(g)}$ vanishes in codimension $g$ it follows
that $[I^{(g)}]=c_g(\EE \otimes \det \EE^{\otimes 1/2} )$. This Chern class can be computed using $c_i(\EE) = \lambda_i$, and then the claim follows, since
there are $2^{g-1}(2^g-1)$ odd theta characteristics.
\end{proof}

\smallskip
Recall that the moduli spaces $\calA_g$ are not compact. By $\Perf$ we denote the perfect cone toroidal compactification of $\calA_g$ and by $\Vor$ we
denote the second Voronoi toroidal compactification.
The rest of this paper is devoted to studying the closure of the locus $I^{(g)}$ in $\Perf$. The boundary degenerations are much harder, and the final result is the computation of the class of the closure in theorem \ref{class}.

\section{Extension of theta gradients on the boundary}
\label{sec:extend}
The boundary of the perfect cone compactification $\Perf$ is an
irreducible divisor $D\subset\Perf$.
We denote by $p:\calA_g(\ell)\to\calA_g$ and $\bar p:\Perfl\to\Perf$ the level $\ell$ covers of the moduli spaces,
and by $D_i$ the irreducible divisorial
components of the boundary of $\Perfl$.  We denote $\delta_i$ the class of the boundary divisor in $CH^1(\Perfl)$.
Note that the cover $\bar p:\Perfl\to\Perf$ branches to order $\ell$ along each $D_i$, and thus we have
${\bar p}_*(\delta_i)=\delta/\ell$.

We have seen that the gradient of the theta function can be interpreted as a section of the vector bundle $\EE\otimes\det\EE^{\otimes 1/2}$. The aim
of this section is to prove that we can extend this to the perfect cone toroidal compactification $\Perf$ of $\calA_g$.

We first remark that the Hodge bundle extends as a vector bundle over any toroidal compactification $\calA_g^{\tor}$ of $\calA_g$ (see \cite{mumhirz}),
as well as over any toroidal compactification of any level cover.
Indeed, to define the Hodge bundle we note that any toroidal compactification admits a universal family of (non-compact)
semiabelian varieties $\calG_g^{\tor}$. This has the ``zero'' section (which is really $1\in\CC^*$ on each torus)
$s: \calA_g^{\tor} \to  \calG_g^{\tor}$, and the Hodge bundle is defined by $\EE:=s^*(\Omega^1_{\calG_g^{\tor}/\calA_g^{\tor}})$.
In particular, the fiber of $\EE$ over $[A] \in \calA_g$ is given by $\EE_{[A]}=\Omega^1_{A,0}$.

In this section we shall work over $\calA_g(4,8)$ and its perfect cone compactification $\Perf(4,8)$.
We will, by abuse of notation, denote the Hodge bundle on this level cover, as well as its extension to the compactification, by $\EE$.

\begin{prop}\label{prop:extension}
The gradients of theta functions with characteristics at zero, $F_m$, extend on $\Perf(4,8)$ to sections of the extension of the Hodge bundle
twisted by a square root of its determinant, i.e.~we have some extensions
$$
 \overline{F_m} \in H^0\left(\Perf(4,8),\det\EE^{\otimes 1/2}\otimes\EE\right).
$$
\end{prop}
\begin{proof}
Before entering into the necessary computations, we will make some comments.
It is enough to prove extension to the generic point on each boundary component, as extension over codimension $2$ sets is then
automatic by Hartogs' extension theorem on normal analytic spaces.
Here we shall make use of the fact that the perfect cone decomposition $\Perf$ has only one boundary component. This is no longer true
for $\Perf(4,8)$, but the group $\Gamma_g / \Gamma_g(4,8)$ acts transitively on the set of boundary components $D_i$ of
$\Perf(4,8)$. Hence it will be sufficient
to consider one of them --- we shall work with the so-called standard boundary component.
Since the Voronoi and the perfect cone compactification coincide in genus $g \leq 3$ and since $\Vor$ is a blow-up
of $\Perf$ in genus $g=4,5$ (\cite{erry2,ryba}) we also obtain extension to the Voronoi compactification for
genus $g\leq 5$.

Recall that the boundary components of $\Perf(4,8)$ correspond to lines in $\QQ^{2g}$ modulo the action of the group
$\Gamma_g(4,8)$. We shall work with the standard cusp corresponding to the line $l_0$
generated by the vector $(0, \ldots , 0; 1, 0, \ldots ,0)$.
Let $P(l_0)$ be the corresponding parabolic subgroup
and let $P'(l_0)$ be the center of the unipotent radical of $P(l_0)$. Moreover let $U(l_0)=P'(l_0) \backslash \calH_g$
be the partial quotient with respect to $P'(l_0)$ and let $V(l_0)$ be the partial compactification of $U(l_0)$ given
by adding the cusp corresponding to $l_0$ (see below for details). The partial compactification of
$\Perf(4,8)$ in a neighborhood of the standard cusp is then obtained by taking the quotient of $V(l_0)$ by the group
$P(l_0)/P'(l_0)$.
Clearly it makes sense to speak about the Hodge bundle $\EE_{\calH_g}$ over
the Siegel space $\calH_g$ (where we have a universal family) as well as about the Hodge bundle $\EE_{U(l_0)}$ resp.
its extension $\EE_{V(l_0)}$ over $U(l_0)$ and $V(l_0)$ respectively. The Hodge bundle $\EE_{\calH_g}$ is trivial as follows
immediately from the construction of the universal family over $\calH_g$. More precisely the universal family over $\calH_g$
is given as the quotient of $\calH_g \times \CC^g$ by the group $\ZZ^{2g}$ where $(M,N) \in \ZZ^{2g}$ acts on $\calH_g \times \CC^g$
by $(\tau,z) \mapsto (\tau, z + \tau M + N)$. The differentials $dz_1, \ldots , dz_g$ then define a trivialization of
$\EE_{\calH_g}$.

\begin{lm} \label{lem:fundamental}
The following holds
\begin{itemize}
\item[(i)] The trivialization of $\EE_{\calH_g}$ over $\calH_g$ descends to a trivialization of $\EE_{U(l_0)}$ over $U(l_0)$.
\item[(ii)] $j_*(\EE_{U(l_0)})=\EE_{V(l_0)}$ where $j: U(l_0) \hookrightarrow V(l_0)$ is the inclusion.
\end{itemize}

\end{lm}
\begin{proof}[Proof of the lemma]
We recall that the universal abelian variety over
$\calA_g(4,8)$ is defined by taking the quotient of $\calH_g \times \CC^g$ with respect to the semi-direct
product $\ZZ^{2g} \rtimes \Gamma_g(4,8)$, which acts as follows:
$$
\left((M,N),\begin{pmatrix} A & B\\C & D \end{pmatrix}\right): \calH_g \times \CC^g \to \calH_g \times \CC^g
$$
$$
(\tau,z) \mapsto ((A\tau + B)(C\tau + D)^{-1}, (z+\tau M + N)(C\tau + D)^{-1}).
$$
The center of the unipotent radical of the parabolic subgroup associated to a line in $\QQ^{2g}$
is a rank $1$ lattice. For the standard cusp
it consists of the matrices $\begin{pmatrix} A & B\\C & D \end{pmatrix}$ of the form
$$
\begin{pmatrix}
1 & 0 & s & 0\\
0 & 1_{g-1} & 0 & 0 \\
0 & 0 & 1 & 0 \\
0 & 0 & 0 & 1_{g-1}\\
\end{pmatrix}
\mbox{ where } s\in 8\ZZ.
$$
In particular $A=D=1_g$ and $C=0$. {}From this we can immediately deduce the first assertion of the lemma.

In order to prove the second assertion we have
to understand how the universal semi-abelian variety can be extended over the generic point of the boundary divisor
associated to the standard cusp. The quotient of $\calH_g$ by $P'$ is given by
$$
\calH_g \to \calH_{g-1} \times \CC^{g-1} \times \CC^*
$$
$$
\tau = \begin{pmatrix} \omega &  b^t\\ b & \tau' \end{pmatrix} \mapsto (\tau', b, \ee(\omega / 8)).
$$
We denote the image of this quotient map by $U$.
This is an open subset (in the analytic topology) of $\calH_{g-1} \times \CC^{g-1} \times \CC$.
Let $V$ be the interior of the closure of $U$ in $\calH_{g-1} \times \CC^{g-1} \times \CC$. The difference between $V$ and $U$ is the
set $\calH_{g-1} \times \CC^{g-1} \times \{0\}$. Adding this set is adding the divisor associated to the standard cusp. We consider
$q_8:=\ee(\omega/8)$ as the coordinate on $\CC^*$. Thus adding the boundary divisor corresponds to adding the set $\lbrace q_8=0\rbrace$.
We now have to understand the universal semi-abelian variety over $U$ and its extension to $V$.
Let $N_1 \cong \ZZ^{2g}$ be the lattice given by $\{M\tau + N\mid M,N \in \ZZ^g\} $ and let
$N_2=\{M; M \in \ZZ^g\} \cong \ZZ^g$.

In order to construct the universal semi-abelian variety we  consider the action
of $N_1 \rtimes P'$ on $\calH_g \times \CC^g$.
For this we first consider the action
of the subgroup $N_2 \rtimes P'$. Clearly, the quotient of $\calH_g \times \CC^g$ by this
subgroup is the trivial $(\CC^*)^g$-bundle on
$U$. We denote the coordinates on $(\CC^*)^g$ by $x_i=\ee(z_i); i= 1, \ldots ,g$. The action
of $N_3=N_1/N_2 \cong {\ZZ}^g$ on $U \times (\CC^*)^g$ is trivial on the base and by
multiplication with powers of $(t_{k,1},t_{k,2}, \ldots ,t_{k,g})$ where $t_{k,j}=\ee(\tau_{kj})$ on the torus $(\CC^*)^g$.
In order to construct the semi-universal abelian variety we extend the trivial $(\CC^*)^g$-bundle on $U$
trivially to $V$ and also extend
the action of $N_3$ to $V \times (\CC^*)^g$.
Note that the action of $N_3$ on $V \times (\CC^*)^g$ is no longer free on the first coordinate when $k=1$.
In order to overcome this difficulty one considers
a toroidal embedding $ V \times (\CC^*)^g \hookrightarrow X_{\Sigma}$. This construction is analogous to the construction of Shioda modular surfaces (see \cite[Part I, 2B and 3D]{hukawebook}) where the case $g=2$ is treated in detail. There is a projection $X_{\Sigma} \to V$ whose restriction to $U$ is just
projection of the trivial torus bundle $U \times (\CC^*)^g \to U$. Over the boundary, i.e.~over points with $q_8=0$ one has a chain of countably many copies
of $(\CC^*)^g$. The action of $N_3$ on $U \times (\CC^*)^g \to U$ extends to an action on $X_{\Sigma}$ and the quotient is a semi-abelian
group scheme over $V$. The semi-abelian group scheme over the partial compactification of $\calA_g(4,8)$ in the direction of the standard
cusp is then obtained by taking a further quotient
with respect to $P/P'$. Note that this group acts freely due to the presence of the level structure.

To prove the second assertion of the lemma we have to understand the relation of $\EE_V$ and  $\EE_U$.
We claim that we can take $dx_1, \ldots ,dx_g$ as a basis for the fibers of $\EE_V$. This indeed proves assertion (ii).
To see the claim it is enough to consider the case of the universal elliptic curve with fiber
coordinate $x_1$, as the other coordinates
are not affected by the construction. Thus we consider the torus $(\CC^*)^2$
with coordinates $(w_1, t_{1,1})=(x,q)$. For this we define
a torus embedding $(\CC^*)^2 \hookrightarrow X'_{\Sigma}$. As we are only interested in the situation over the section given by
the origin it suffices to look at one chart $X'_{\sigma_0}$ of $X'_{\Sigma}$. We are thus in the situation of
\cite[p.~29]{hukawebook} and $X'_{\sigma_0} \cong \CC^2$ with
embedding $(x,q) \mapsto (x,x^{-1}q)=(u,v)$. The projection onto the
base is given by $(u,v) \mapsto uv=q$. The section given by the origin is the set $\{(1,q); q\in \CC\}$. In order to describe the Hodge bundle
we have to consider the relative tangent bundle $\Omega^1_{\CC^2/\CC}$ restricted to the zero-section. This is generated by $du, dv$ modulo
the pullback $d(uv)=udv + vdu$. Restricting the latter to the zero-section
gives $dv + vdu$ and hence the relative cotangent bundle
is generated by $du=dx$ which we can take as trivializing section.
\end{proof}

This argument shows that, in particular, the Hodge bundle over $V$ is trivial.
Hence all we have to do to prove our claim is to show that the functions  $F_m$ are invariant with respect to $P'$ (which is obvious) and that
they extend without poles to $V$. For this we have to compute the Fourier-Jacobi expansion. This was done in detail in \cite{grsmconjectures}, and we now summarize the results for completeness. Indeed, denote $z=(z_1,z')$ where $z_1 \in \CC$ and $z' \in \CC^{g-1}$ and write $
\tau = \left(\begin{matrix} \omega& b^t \\ b &\tau' \end{matrix}\right).
$
The Fourier-Jacobi expansion at the standard cusp amounts to writing the Taylor series in $q_8$ (the coordinates transverse to the boundary of $\Perf(8)$) as $\tau_{11}\to i\infty$ For the characteristics of the gradient being $\ep=\ep_1 \ep'$ and $\de=\de_1 \de'$
we then get for the case $\ep_1=0$
\begin{equation}\label{equ:vanishing1}
  \partial_{z_1}\tc{0&\ep'}{\de_1&\de'}=q_8^4\cdot 4\pi i  (-1)^{\de_1}\tc{\ep'}{\de'}(\tau,b)+O(q_8^{16})
\end{equation}
\begin{equation}\label{equ:vanishing2}
  \partial_{z}\tc{0&\ep'}{\de_1&\de'}=2\partial_{z}\tc{\ep'}{\de'}(\tau,z)|_{z=0}+O(q_8^4).
\end{equation}
This is to say that in the case of $\ep_1=0$ (this can be said invariantly: see Proposition \ref{propF} below)
the gradient does not vanish identically. On the other hand, for $\ep_1=1$ we get
\begin{equation}\label{equ:vanishing3}
  \partial_{z_1}\tc{1&\ep'}{\de_1&\de'}=q_8\cdot 4\pi i  \ee(\de_1/4)\tc{\ep'}{\de'}(\tau,b/2)+O(q_8^9)
\end{equation}
\begin{equation}\label{equ:vanishing4}
  \partial_{z}\tc{1&\ep'}{\de_1&\de'}=q_8\cdot 2\ee(\de_1/4)\partial_{z}\tc{\ep'}{\de'}(\tau,b/2)+O(q_8^9)
\end{equation}
which shows that in this case the generic vanishing order of $\overline{F_m}$ in $q_8$ is precisely equal to one. In any case the Fourier-Jacobi expansion is
holomorphic and in view of Lemma \ref{lem:fundamental} this proves Proposition \ref{prop:extension}.
\end{proof}

The above computations actually give us more information.
To explain this, we prefer to work on the full level-$8$ cover $\Perf(8)$ rather than $\Perf(4,8)$. Recall from \cite[Ch.~4]{namikawabook} (or eg.~\cite[Sec.~3]{erdenbergerthesis}) that the
boundary components of $\Perf(8)$ are in bijective correspondence
with the primitive vectors in $((\ZZ/8\ZZ)^{2g} \setminus \{0\})/\pm 1$, and $\Gamma_g$ acts transitively on the set of boundary components of $\Perf(8)$.
Under this correspondence the standard boundary component corresponds to the vector $(0,0,\ldots,0;1,0,\ldots,0)$.
Now let $D_n$ be a boundary component of $\Perf(8)$ labeled by
a primitive vector $\pm n \in (\ZZ/8\ZZ)^{2g}$. We  denote the reduction of $n$ modulo 2 (which simply sends odd entries to
$1$ and even entries to $0$) by $n_2$.
We identify an odd 2-torsion point $m=(\tau\ep +\de)/2$ with its characteristic $(\ep,\de)\in (\ZZ/2\ZZ)^{2g}$.
The above calculations show that $F_m$ does not vanish on the standard component if and only if $\ep_1=0$.
Note that in this case $\sigma(m+n_2)=\ep\de+ \ep_1=1$ since $m$ is odd, i.e. $\ep\de=1 \mod 2$ and $\ep_1=0$.
Using the action of the finite group $\Gamma_g/\Gamma_g(8)$ on $\Perf(8)$
we thus obtain

\begin{prop}\label{propF}
The gradient $\overline{F_m}$ does not vanish generically on the
boundary component $D_n$ of $\Perf(8)$ if and only if $\sigma(m+n_2)=1$.
Otherwise it vanishes on $D_n$ with generic vanishing order equal to one.
\end{prop}
We would like to interpret the above proposition by saying that
\begin{equation}\label{tildeF}
 \tilde F_m\in H^0\left(\Perf(8),\det(\EE)^{\otimes1/2}\otimes\EE\otimes \calO\left(-\sum\limits_{\lbrace \pm n\mid m+n_2\rm \even\rbrace}D_n\right)\right).
\end{equation}
To be able to say this we must show that $\sum\limits_{\lbrace \pm n\mid m+n_2\rm \even\rbrace}D_n$ is a Cartier divisor
and not just a Weil divisor. This is not immediate as $\Perf(8)$ is not a smooth space for $g \geq 4$.
Analytically this means that we want to divide the gradient by the product of the defining
equations for all boundary components
of $\Perf(8)$ where it
vanishes identically with multiplicity one.

The gradients of theta functions with characteristics were studied in \cite{grsmodd1}, and their boundary
degenerations were also considered in \cite{grsmconjectures}. It seems, however, that
the gradients $f_m$ of the theta function at odd two-torsion points were not
as extensively studied, while for us they turn out to be more important.

\begin{prop}\label{Fandf}
The divisor $\sum\limits_{\lbrace \pm n\mid m+n_2\rm \even\rbrace}D_n$ is a Cartier divisor
and we can also extend the sections $f_m$ to sections
$$
 \tilde f_m\in H^0\left(\Perf(8),\det(\EE)^{\otimes(1/2)}\otimes\EE\otimes \calO\left(
-\sum\limits_{\lbrace \pm n| m+n_2\rm \even\rbrace}D_n\right)\right)
$$
not vanishing identically on any component $D_i$ of the boundary of $\Perf(8)$.
Moreover $\tilde f_m$ and $\tilde F_m$ are equal up to a non-zero constant and thus we have
$$
 \overline{G_m^{(g)}}(8):=\lbrace \tau\in\Perf(8)\mid \tilde F_m=0\rbrace
 =\lbrace \tau\in\Perf(8)\mid \tilde f_m=0\rbrace.
$$
\end{prop}
\begin{proof}
Indeed, by formula (\ref{fmdef}) we have on $\calA_g(8)$ the relation
$$
 f_m(\tau)=\ee(-\ep^t\tau\ep/8-\ep^t\de/4)F_m(\tau),
$$
and thus to
determine the extension of $f_m$ to $\Perf(8)$ it is enough to determine
the properties of the exponential factor, and to combine them with
the results obtained above about $\overline{F_m}$.
Note also that the functions $\ee(-\ep^t\tau\ep/8-\ep^t\de/4)$ descend to any partial quotient which arises in
the construction of the toroidal compactification and thus $\sum\limits_{\lbrace \pm n\mid m+n_2\rm \even\rbrace}D_n$
is a Cartier divisor on $\Perf(8)$.
To determine the extension of $f_m$  it is again enough to work with the
standard component. Noticing
that $\ee(-\ep^t\de/4)$ is a constant independent of $\tau$ and thus irrelevant
for our computations, we compute in the above notation
$$
 \ee\left(-\frac18\begin{pmatrix} \ep_1& {\ep'}^t\end{pmatrix}
 \begin{pmatrix} \omega &  b^t\\ b & \tau' \end{pmatrix}
 \begin{pmatrix} \ep_1\\ \ep'\end{pmatrix}\right)
$$
$$
 =\ee\left(-\frac18\left(\ep_1^2\omega+2\ep_1b^t\ep'+{\ep'}^t\tau'\ep'\right)\right)=
 q_8^{-\ep_1^2}\ee(-\ep_1b^t\ep'/4-{\ep'}^t\tau'\ep'/8).
$$
It thus follows that for $\ep_1=0$ this exponential factor is independent of $q_8$,
while for $\ep_1=1$ it has a pole $q_8^{-1}$. Noticing that in these
cases $\overline{F_m}=O(1)$ and $\overline{F_m}=O(q_8)$, respectively, it thus follows that the extension
$\tilde f_m$ does not vanish generically on any boundary component of $\Perf(8)$, and thus finally that
$\tilde f_m$ and $\tilde F_m$ coincide up to a non-zero constant.
\end{proof}
Since the group $\Gamma_g/\Gamma_g(8)$ permutes the sets $\overline{G^{(g)}_m(8)}\subset\Perf(8) $, this defines a locus
$\overline{G} =\overline{G^{(g)}}$ such that $\cup_m \overline{G^{(g)}_m(8)}=
p^*(\overline{G^{(g)}})$.

Geometrically, the above computation amounts to saying that we have determined the vanishing orders of $\overline{F_m}$ on the open part of the boundary of $\Perf(8)$, and in particular noticed that the zero loci of $\tilde{f_m}$ do not have any irreducible components contained in this locus. It is then natural to wonder whether if we consider the full compactification $\Perf(8)$, the zero locus of any $\tilde{f_m}$ may have any irreducible component contained in the boundary (from the above, it would then have to be disjoint from the partial compactification). We conjecture that this in fact impossible:
\begin{conj}\label{conj:setG}
For any genus we have the equality
$
\overline{G^{(g)}}= \overline{I^{(g)}}.
$
\end{conj}
While we were unable to prove this conjecture in full generality, one of the main results of \cite{grhu2} is theorem 1.3 there, which is a proof of the above conjecture for $g\le 5$. This result is obtained by explicitly describing the geometry of various types of principally polarized semi-abelic varieties in detail, and checking that there are no components of $\overline{G^{(g)}}$ contained in the loci of semi-abelic varieties of any of the types described. We thus obtain our second main result:
\begin{proof}[Proof of theorem \ref{class}]
Recall that what this theorem claims is an expression for the class of the locus $\overline{I^{(g)}}$ in $CH^g(\Perf)$. Note, however, that the locus $\overline{G^{(g)}_m(8)}\subset\Perf(8)$ is by definition the zero locus of $\tilde{f_m}$, which is a section of a certain vector bundle given by proposition \ref{Fandf}. Thus the class of $\overline{G^{(g)}_m(8)}$, if it has expected codimension $g$, is given by the top Chern class of this vector bundle. Thus we need to first understand the geometry of the cover, and then to compute the top Chern class.

Indeed, we are working with the full level $\ell$ Galois cover $ p: \Perfl \to \Perf$, with the Galois group $\Sp(2g,\ZZ/\ell\ZZ)$. As long as one considers projective varieties this action factors through
$\PSp(2g,\ZZ/\ell\ZZ)$.
However, as we will be doing all our calculations in the stack set-up it is important to
work with $\Sp(2g,\ZZ/\ell\ZZ)$ as $\pm {\bf 1}$ acts non-trivially on $\Perf$ (as a stack).
For any $\Sp(2g,\ZZ/\ell\ZZ)$-invariant cycle $Z$ with class $[Z]$ on $\Perf (n)$ we define
\begin{equation}\label{projection}
{\overline p}_*([Z]) = \frac{1}{|\Sp(2g,\ZZ/\ell\ZZ)|} [p_*(Z)].
\end{equation}
For classes $Z$ which are the pullback of a stacky class $W$ we have  ${\overline p}_*(Z)=W$ and in particular
we have ${\overline p}_*(\lambda_i)= \lambda_i$ where the $\lambda_i$ are the classes of the Hodge bundle(s).

Now, to prove the theorem we first perform the computation on $\Perf(8)$, and then argue invariance to reduce this to $\Perf(2)$. Denote by $\al_1,\cdots,\al_g$ the Chern roots of the Hodge vector bundle $\EE$ on $\Perf(8)$, so that the
Chern polynomial is $c(\EE)=\prod (1+\al_i)$, which means that for the symmetric polynomials we have
$s_i(\al_1,\ldots,\al_g)=\lambda_i$. {}From formula (\ref{tildeF}) we know that $\tilde f_m$ is a section
of $\EE\otimes L_m$ for a certain line bundle $L_m$.
The cycle $\overline{G^{(g)}}(8)$ (which is the pre-image of  $\overline{G^{(g)}}$)
is the union of the vanishing loci of the sections $F_m$.
Now, if $\tilde f_m$ vanish in codimension $g$, so that $\overline{G}$ is pure of codimension $g$ --- and this is our assumption (or the content of \cite[Thm.~1.2 and 1.3]{grhu2} for $g\le 5$), we recall that the vanishing locus is reduced, as noted in remark \ref{remred}.
It now follows, see \cite[Chapter 14]{fultonintersection}, that
$$
 [\overline{G^{(g)}}(8)]=\sum\limits_{m\in(\ZZ/2\ZZ)^{2g}_{\rm odd}} c_g(\EE\otimes L_m).
$$
Denoting
$$
 \ell_m:=c_1(L_m)= \frac{\lambda_1}{2}-\sum\limits_{\pm n\in (\ZZ/8\ZZ)^{2g}_{\rm primitive}, n_2+m\rm\ even}\delta_n
 =\frac{\lambda_1}{2}-\sum\limits_{n\in Z_m}\de_n
$$
we compute on $\Perf(8)$ that
$$
 c_g(\EE\otimes L_m)= \left[\prod\limits_{ i=1}^g(1+\al_i+\ell_m) \right]_{\deg=g}=\sum\limits_{i=0}^g s_{g-i}(\al_1,\ldots,\al_g)\ell_m^{i}
$$
$$
 =\sum\limits_{i=0}^g\lambda_{g-i} \ell_m^i=\sum\limits_{i=0}^{g} \lambda_{g-i}\left(\frac{\lambda_1}{2}-\sum\limits_{n\in (\ZZ/8\ZZ)^{2g}_{\rm primitive}, n_2+m\rm\ even} \de_n\right)^i
$$
and the formula on $\Perf(8)$ follows.

Finally, to obtain the stated formula, as a pushforward from $\Perf(2)$, we note that the expression above is certainly invariant under the Galois action $\Gamma_g(2)/\Gamma_g(8)$, as only the reduction $n_2$ modulo two matters. Remembering that $\Perf(8)\to\Perf(2)$ branches along the boundary to order 4, we get the result as stated.
Note also that as we remarked before, the vanishing locus of each $F_m$ is well-defined as a subvariety
of $\calA_g(2)$, even though $F_m$ itself is only a section of a well-defined vector bundle over $\calA_g(4,8)$.
\end{proof}

\section{The tautological ring}
We will now use various properties of the tautological ring of $\Perf$ (this section uses the perfect pairing, and we thus have to work with rational coefficients) to compute the projection of the class $[\overline{I^{(g)}}]$ given by theorem \ref{class} to the tautological ring, thus proving theorem \ref{thmtaut}.

We recall that the perfect cone (and in fact any toroidal) compactification $\Perf$ admits a stratification induced by the natural map $P:\Perf\to\Sat$ to the Satake compactification. Writing the Satake compactification as
$$
 \Sat = \calA_g \sqcup \calA_{g-1} \sqcup \calA_{g-2} \ldots \sqcup \calA_0,
$$
following van der Geer \cite{vdgeerchowa3}, we denote $\beta_k^0:= P^{-1}(\calA_{g-k})$ the open strata, and denote by $\beta_k:=\beta_k^0\sqcup\ldots\sqcup\beta_g^0=P^{-1}(\calA_{g-k}^{\rm Sat})$ the closed strata --- and by abuse of notation, their classes in $CH^k(\Perf)$. We also recall that $\delta_i$ denote the classes of the irreducible components of the boundary of $\Perf(2)$, and we follow van der Geer in denoting by $\sigma_k$ the degree $k$ elementary symmetric polynomial in the classes $\delta_i$, which we would like to
consider as an element of $CH^k(\Perf(2))$. This, however, is problematic, as in general the boundary components
$\delta_i$ are not Cartier divisors (this can be seen for $g=4$ by looking at the second perfect cone which is a $10$-dimensional cone with $12$ extremal rays, see \cite{husaA4} for a discussion of this), and thus this is also the case for $g>4$. However, the sum of all boundary components, as well as certain partial sums, are Cartier as we saw in section \ref{sec:extend}.
In what follow, whenever we work with expressions such as the $\sigma_i$ we shall mean by this the restriction of this cycle
to the smooth part $\calA_g^{\op{perf},0}$ of $\Perf$ where it is well defined. We observe that the singular locus of $\Perf$ is
contained in $\beta_4$.
We also notice that $\sigma_k$ is invariant under the action of $\Gamma_g$, and we can thus also denote $\sigma_k$ the corresponding class in $CH^k(\calA_g^{\op{perf},0})$.

We shall frequently make implicit use of the following very useful fact, see \cite[Proposition 1.1.8]{fultonintersection}.
If $Z$ is a closed subvariety of $\Perf$ then there is an exact sequence
\begin{equation}\label{equ:restriction}
CH_k(Z) \to  CH_k(\Perf) \to CH_k(\Perf - Z) \to 0.
\end{equation}

\smallskip
We start with two useful lemmas concerning the projection of classes supported on the boundary to the tautological ring.
\begin{lm}\label{projtotaut}
Let $[X]$ be a class of codimension at most $g-1$ in
the Chow group of a toroidal compactification $\calA_g^{\op{tor}}$, supported over the boundary $\beta_1$.
Then $(P(\lambda)[X])^{\op{taut}}=0$ for every non-constant polynomial $P(\lambda)$ in the tautological classes $\lambda_i$.
\end{lm}
\begin{proof}
We first notice that it is enough to prove that $[X]^{\op{taut}}=0$: if $[X]$ pairs to zero with every polynomial of
complementary degree in the $\lambda_i$, then the same is true for $P(\lambda)[X]$.

We write $[X]= Q(\lambda) + [Y]$ where $Q(\lambda)$ is a polynomial in the $\lambda_i; i\leq g-1$ by the assumption
on the codimension of $[X]$, and $[Y]$ pairs to zero with all tautological classes of complementary degree (that is, polynomials in the $\lambda_i$ of complementary degree).
Restricting to $\calA_g$ we obtain $[X]|_{\calA_g}=0$ and $[Y]|_{\calA_g}$ still projects to $0$ in the tautological ring
of $\calA_g$. But this is a contradiction, since $Q(\lambda)$ defines a non-zero class in
the tautological ring of $\calA_g$, the only new relation in the tautological ring of $\calA_g$ versus that in the tautological ring of $\calA_g^{\op{tor}}$ being $\lambda_g=0$ in degree $g$ (see \cite{vdgeercycles}).
\end{proof}

\begin{lm}\label{projzero}
Let $[X]$ be a class of codimension at most $g$ in the Chow group of a toroidal compactification $\calA_g^{\op{tor}}$ supported over $\beta_2$. Then $[X]^{\op{taut}}=0$.
\end{lm}
\begin{proof}
The idea of the proof is similar to the proof of the previous lemma --- basically these lemmas together are based on the fact that in a sense $\lambda_g$ is the only degree $g$ tautological class supported on $\beta_1$. Rigorously, suppose that we have $[X]=Q(\lambda)+[Y]$ in $CH^*(\Perf)$ with $[Y]^{\op{taut}}=0$. Unless we have $Q(\lambda)=c\lambda_g$, for some $c\in\CC$, the image of $Q(\lambda)$ in $CH^*(\calA_g)$ (obtained from $CH^*(\Perf)$ by imposing one extra condition $\lambda_g=0$) is non-zero, which contradicts the fact that the class is supported over $\beta_2$. But if we have $Q(\lambda)=c\lambda_g$, then we must have $$0=\lambda_1^{\frac{g(g-1)}{2}}[X]=c\lambda_1^{\frac{g(g-1)}{2}}\lambda_g+0\ne 0$$
where we have used the fact that $\lambda_1^{\frac{g(g-1)}{2}}$ is zero on $\beta_2$ for dimension reason, and the last intersection number is non-zero by Hirzebruch-Mumford proportionality. (This can also be seen explicitly without a reference to proportionality, since by \cite[Theorem 4.7]{vdgeercycles} the class $\lambda_g$ is proportional to the class $[B_g]$ of the closure of the zero section of $\calX_{g-1}\to\calA_{g-1}$, modulo classes supported on $\beta_2$, which have zero pairing with the top power of $\lambda_1$ for dimension reasons. By \cite{shepherdbarron} the normalization of this closure of the locus of trivial extensions is isomorphic to $\calA_{g-1}^{\op{Perf}}$, where the top self-intersection number of $\lambda_1$ is non-zero).
\end{proof}

We now compute the projection to the tautological ring of the one degree $g$ monomial in $\sigma$ and $\lambda$ not covered by the above lemmas.
\begin{prop}
The projection of $\sigma_1^g$ to the tautological ring is
\begin{equation}\label{Dg}
[\sigma_1^g]^{\op{taut}}=(-2)^{g-1}(g-1)![B_g]^{\op{taut}}=\frac{-2^{g-1}(g-1)!}{\zeta(1-2g)}\lambda_g,
\end{equation}
where $B_g$ is the class of the closure of the zero section of $\calX_{g-1}\to\calA_{g-1}$ in the boundary of the partial compactification, as above.
\end{prop}
\begin{proof}
We first observe that it suffices to prove the first equality, since it follows from
\cite[Theorem 4.7]{vdgeercycles} and the above lemma \ref{projzero} that $[B_g]^{\op{taut}}=(-1)^g\lambda_g/\zeta(1-2g)$.
Recall $\sigma_1|_{\beta_1^0}=-2T$, where $\beta_1^0=\beta_1\setminus\beta_2$ is the open part (isomorphic to the universal Kummer family),
and $T\subset\calX_{g-1}$ is the universal theta divisor trivialized along the zero section (see \cite{ergrhu2} for a more precise description).
Thus we compute $\sigma_1^g=\overline{(-2T)^{g-1}}+X$ for some class $X$ supported within $\beta_2$.
Here $\overline{(-2T)^{g-1}}$ is an extension of the class $(-2T)^{g-1}$ to $\Perf$ (which is unique up to
a cycle supported on $\beta_2$).
Since the class $[\sigma_1^g] - (-2)^{g-1}(g-1)![B_g]$ is supported on $\beta_1$ we can argue as in the proof of
lemma \ref{projzero} that $[\sigma_1^g]^{\op{taut}} - (-2)^{g-1}(g-1)![B_g]^{\op{taut}}=c \lambda_g$.
Now note that on each fiber of the boundary of the partial compactification $\calX_{g-1}\to\calA_{g-1}$ the divisor $T$ restricts to the principal polarization, with top self-intersection number $(g-1)!$. Thus the degree of $T^{g-1}$ on each fibre is $(g-1)!$. Intersecting with $\lambda_1^{\frac{g(g-1)}{2}}$ then shows that $c=0$ as claimed.
\end{proof}

We are now ready to compute the projection $[\overline{I^{(g)}}]^{\op{taut}}$.
\begin{proof}[Proof of theorem \ref{thmtaut}]
Let us expand the expression in formula (\ref{Gclass}) in theorem \ref{class}, and take one of the terms: this will then be of the form $p_*(Q(\lambda)P(\delta_n))$, where $Q(\lambda)$ is a monomial in the classes $\lambda_i$ (i.e.~lies in the tautological ring), $P$ is a monomial in classes $\delta_n$ of various irreducible components of the boundary of $\Perf(2)$, and the total codimension is equal to $g$.

If $P$ contains a product of at least two different boundary components, $\de_{n_1}\de_{n_2}$, it is supported within $\beta_2$, and thus by lemma \ref{projzero} in this case we have $[P]^{\op{taut}}=0$. This implies that all pairing of $P$ with polynomials in $\lambda$ classes are zero, and thus this is also the case for $QP$, so we have $[p_*(Q(\lambda)P(\de_n))]^{\op{taut}}=0$. Thus to have a non-zero projection to the tautological ring we must have $P=\de_n^i$. Moreover, if $0<i<g$, then lemma \ref{projtotaut} applies, and we also get $[p_*(Q(\lambda)P(\de_n))]^{\op{taut}}=0$.
Thus the only monomials in the expansion of (\ref{Gclass}) that can have a non-zero projection to the tautological ring are either those that do not contain any $\de_n$ (and thus these are the ones giving the expression on the open part in theorem \ref{thm:Igclass}), or the ones of the form $\de_n^g$. Thus the extra term compared to theorem \ref{thm:Igclass} that we need to compute is
$$
  \frac{1}{N}\sum\limits_{m\in(\ZZ/2\ZZ)^{2g}_{\rm odd}}(-4)^{-g}{p}_* \left(\sum\limits_{n\in Z_m} \de_n^g\right).
$$
We know by proposition \ref{propF} that $n\in Z_m$ if and only if $\sigma(n+m)=0$. Thus for $n$ fixed there are $2^{2g-2}$ such odd $m$.
The projection $p_*\de_n$ for any $n$ is equal to $\sigma_1/2$ (recall the branching order 2 of the cover $\Perf(2)\to\Perf$ along the boundary). Thus the above pushforward is equal to
$$
 \frac{1}{N}(-1)^g2^{-2g}2^{2g-2}\sum\limits_{n\in(\ZZ/2\ZZ)^{2g}\setminus\{0\}}p_*(\de_n^g)
 =(-1)^g2^{-2-g}\sigma_1^g.
$$
We now use formula (\ref{Dg}) to obtain
$$
(-1)^g2^{-2-g}[\sigma_1^g]^{\op{taut}}= \frac{(-1)^{g-1}(g-1)!}{8\zeta(1-2g)}\lambda_g+[X]^{\op{taut}},
$$
where, however, the class $[X]$ is of codimension $g$, supported over $\beta_2$, and thus has zero tautological projection by lemma \ref{projzero}.
\end{proof}

\section{Combinatorics of boundary divisors}
The rest of this paper will be devoted to using theorem \ref{class} to compute the classes of the loci $\overline{I^{(g)}}$ explicitly for $g\le 5$ (and not only its tautological part given by theorem \ref{thmtaut}).

In this section we will deal with the combinatorics of the terms appearing in the expression in theorem \ref{class}, enabling us to further perform the explicit class computations.
\begin{prop}
Formula (\ref{Gclass}) can be written in a different way, grouping by powers of $\delta$, as follows:
\begin{equation}\label{Gclass2}
 [\overline{I^{(g)}}]=\sum\limits_{m{\rm\ odd}}{\overline p}_*\left(\sum\limits_{j=0}^g \left(
 -{\sum\limits_{n\in Z_m} \de_n}\right)^j\sum\limits_{i=j}^g
\binom{i}{j}\left(\frac{\lambda_1}{2}\right)^{i-j}\lambda_{g-i}\right).
\end{equation}
\end{prop}
This expression is easier to use, as similar terms involving $\delta_n$'s are
grouped together.
To use the proposition we thus need to study the combinatorics of the intersections of the components $D_i$
of the boundary
$\partial\Perf(2)$, and the vanishing of a given $\overline{F_m}$ on a set of components.

The combinatorial description of the boundary components of $\Perf(2)$ and their intersection behavior can be deduced
from the correspondence between cusps and isotropic subspaces of $\QQ^{2g}$.
For a detailed discussion of the combinatorics of the boundary we refer the reader to \cite{erdenbergerthesis}.
Recall that the group
$\Sp(2g,\ZZ/2\ZZ)=\Gamma_g/\Gamma_g(2)$ has an affine transitive action on the set of
odd theta characteristics $m$. The action of $\Sp(2g,\ZZ/2\ZZ)$ on the set of
boundary components of $\Perf(2)$, i.e.~the linear
action on $(\ZZ/2\ZZ)^{2g}\setminus\{0\}$ is also transitive.
Recall that the symplectic product of $n_1,n_2\in(\ZZ/2\ZZ)^{2g}$ with $n_i=[\al_i,\be_i]$
is defined as $n_1\cdot n_2=\al_1\cdot\be_2+\al_2\cdot\be_1\in\ZZ/2\ZZ$.
Then boundary components $D_{n_1}$ and $D_{n_2}$ of $\Perf(2)$ intersect if and only
if $n_1\cdot n_2=0$, i.e.~$n_1, n_2$ span an isotropic plane (see. \cite[Proposition 3.3.15]{erdenbergerthesis}).
Moreover recall that the intersection $D_{n_1}\cap\ldots\cap D_{n_k}\subset\partial\Perf(2)$ lies in the
stratum $\beta_j$ and intersects $\beta_j^0$, where $j$ is the dimension of the linear span of the
vectors $n_1,\ldots,n_k$ in $(\ZZ/2\ZZ)^{2g}$ (i.e.~an open part of this intersection lies in the preimage of
$\calA_{g-j}\subset\partial\Sat$).

\begin{prop}\label{intersect}
For any set of vectors $n_1,\ldots,n_k\in(\ZZ/2\ZZ)^{2g}\setminus\{0\}$ such that the intersection
$D_{n_1}\cap\cdots\cap D_{n_k}$ is non-empty
the number of gradients $F_m$ vanishing on each $D_{n_i}$, i.e.~the number of odd
$m\in(\ZZ/2\ZZ)^{2g}$ such that $m+n_i$ is even for all $i$, is
\begin{enumerate}
\item zero if there exists a linear relation $n_{i_1}+\ldots+n_{i_j}=0$ with odd number of terms $j=2l+1$
\item equal to $2^{2g-k-1}$ if $n_1,\ldots, n_k$ are a basis for an isotropic subspace
\end{enumerate}
\end{prop}
\begin{proof}
For part (1), relabel the $n$'s so that the linear dependence is $n_1+\ldots+n_{2l+1}=0$. For the intersection $D_{n_1}\cap\ldots\cap D_{n_{2l+1}}$ to be non-empty, we must have
$n_i\cdot n_j=0$ for all $1\le i<j\le 2l+1$. On the other hand, for $\overline{F_m}$ to vanish on $D_{n_i}$, we must have $n_i\in Z_m$, i.e.~$\sigma(m+n_i)=0$.
Recall that by definition we have $\sigma(m+n)=\sigma(m)+\sigma(n)+m\cdot n$, and thus
$$
 \sigma(m_1+\ldots +m_k)=\sum\limits_{i=1}^k\sigma(m_i)+\sum\limits_{ 1\le i<j\le k} m_i\cdot m_j.
$$
Suppose now that some $\overline{F_m}$ vanished on all of these $D_{n_i}$.
We can then use the linear dependence to compute
$$
  1=\sigma(m)=\sigma\left(m+\sum\limits_{i=1}^{2l+1} n_i\right) =
$$
$$
\sigma(m)+\sum\limits_{i=1}^{2l+1}\sigma(n_i)
  +m\cdot\sum\limits_{i=1}^{2l+1} n_i+\sum\limits_{1\le i<j\le 2l+1}n_i\cdot n_j=
$$
$$
 =(2l+1)\sigma(m)+\sum\limits_{i=1}^{2l+1}\sigma(n_i)+m\cdot\left(\sum\limits_{i=1}^{2l+1} n_i\right)=
 \sum\limits_{i=1}^{2l+1}\sigma(m+n_i)=0,
$$
which is a contradiction.

For part (2), note that the action of $\Sp(2g,\ZZ/2\ZZ)$ is transitive on the set of bases of
$k$-dimensional isotropic subspaces (this is well known, eg. \cite[Remark 3.3.9]{erdenbergerthesis}),
and thus any such basis $k$-tuple is conjugate to the standard one, i.e.~in which $n_i=[0,e_i]$ with $e_i$ being the unit vector in $(\ZZ/2\ZZ)^g$ in the $i$'th direction. Thus it is enough to count the number of $m$ such that $\sigma(m+n_i)=0$ for the standard isotropic tuple. Indeed, such $m$ are all of the form
$$
\left[\begin{matrix} \ep\\ \de \end{matrix} \right] =
 \left[\begin{matrix} 1\ldots 1 *\ldots *\\ {*\ldots*} *\ldots *\end{matrix}\right]
$$
where $*$ denote arbitrary entries and there are $k$ 1's in the top row. Recalling that $m$ must be odd means we can choose all entries except the first bottom one arbitrarily, and then there would be a unique choice of the first bottom one to make $m$ odd. Thus the total number of such $m$ is equal to $2^{2g-k-1}$ as claimed.
\end{proof}

\section{Known cases: genus 2 and 3}\label{sec:lowgenus}
We shall now demonstrate how our formula works by computing the classes $[\overline{G^{(g)}}]$ in genus 2 and 3.  Note that in general it is an extremely interesting question to define a meaningful ``enlarged'' tautological ring of the Chow ring for some toroidal compactification of $\calA_g$ (see \cite{vdgeercycles},\cite{ekvdgcycles}, and we hope our computations here may also shed further light on this.

Note that for $g\le 3$ the stack $\calA_g^{\op{perf}}$ is smooth and thus we can work in the (classical) Chow ring of this space.
\subsection{Genus $2$}
Recall that the locus $I^{(2)}$ is empty, and also $ \overline{G^{(2)}}=\emptyset$.
It is convenient to introduce some ordering on the set $(\ZZ/2\ZZ)^{2g}$, so that the following formulas are easier to write using $n_1<n_2$, etc.~notation.
We now use Theorem \ref{class} to compute
$$
 [\overline{G^{(2)}}]=\sum\limits_{m{\rm\ odd}}{\overline p}_* \left(
 \lambda_2+ \lambda_1\left(\frac{\lambda_1}{2} -\frac14\sum\limits_{n\in Z_m} \de_n\right) + \left( \frac{\lambda_1}{2} -\frac14
\sum\limits_{n\in Z_m} \de_n\right)^2\right)
$$
$$
 ={\overline p}_*\left(\sum\limits_{m\rm\ odd} \left( \frac{5}{4}\lambda_1^2-\frac24\lambda_1 \sum\limits_{n\in Z_m} \de_n+\frac{1}{16}\sum\limits_{n\in Z_m} \de_n^2+\frac{2}{16}\!\!\sum\limits_{n_1<n_2 \in Z_m; n_1\cdot n_2=0}\!\!\!\!\!\!\de_{n_1}\de_{n_2}\right)\right)
$$
where we used the relation $\lambda_2= \lambda_1^2/2$.
In genus 2 there are 6 odd characteristics. By proposition \ref{intersect} four of $\overline{F_m}$ vanish on any boundary component $D_i$. Any intersection $D_i\cap D_j$ of two boundary components, if non-empty, is codimension two, and there exist then two $\overline{F_m}$ vanishing on both $D_i$ and $D_j$. We thus compute
$$
 [\overline{G^{(2)}}]={\overline p}_*\!\left(
6\cdot\frac{5}{4}\lambda_1^2-4 \cdot \frac12\lambda_1\sum_n \de_n+ 4\cdot\frac{1}{16}\cdot\sum_n \de_n^2+2\cdot\frac18\sum_{n_1< n_2}\de_{n_1}\de_{n_2} \right)
$$
$$
 =\frac{15}{2}\lambda_1^2-\lambda_1\sigma_1+\frac{\sigma_1^2- 2
 \sigma_2}{16}+\frac{\sigma_2}{16}
$$
where we again used the fact that $\Perf(2)\to\Perf$ has branching order 2 along the boundary, so that the
pushforwards are
\begin{equation}\label{push1}
  {\overline p}_*(\sum \de_n)=\frac{\sigma_1}{2}, \qquad {\overline p}_*(\sum_{n_1< n_2} \de_{n_1}\de_{n_2})=\frac{\sigma_2}{4}
\end{equation}
as well as
\begin{equation}\label{push2}
 {\overline p}_*\left(\sum_n \delta_n^2\right)={\overline p}_*\left((\sum_n \delta_n)^2-2\sum_{n_1< n_2} \delta_{n_1}\delta_{n_2}\right)=\frac{\sigma_1^2}{4}-\frac{\sigma_2}{2}.
\end{equation}
Here $\sigma_2$ is as in \cite{vdgeercycles}, namely the stacky cycle defined by the union of the intersections
of two boundary components on a level cover.
Finally from the computation of $CH^*(\calA_2^{\op {Perf}})$ in \cite{vdgeerchowa3} we have
$$
 \sigma_2=6\lambda_1\sigma_1, \qquad \sigma_1^2=22\sigma_1\lambda_1-120\lambda_1^2
$$
and substituting this into the expression above gives $[\overline{G^{(2)}}]=0$ as expected.

In general these computations allow us to compute the terms with $i<3$ in (\ref{Gclass2}), obtaining for any genus
\begin{lm}\label{classlm}
For any genus we have the formula
$$
 [\overline{G^{(g)}}]= 2^{g-1}(2^g-1)\sum\limits_{i=0}^g\lambda_{g-i} \left(\frac{\lambda_1}{2}\right)^i
 -2^{2g-5}\sigma_1\sum\limits_{i=1}^g i\lambda_{g-i}\left(\frac{\lambda_1}{2}\right)^{i-1}
$$
$$
 +2^{2g-8}(\sigma_1^2-\sigma_2)\sum\limits_{i=2}^g\frac{i(i-1)}{2}\lambda_{g-i}
 \left(\frac{\lambda_1}{2}\right)^{i-2}+P
$$
where the last term denotes some polynomial in $\delta_n$ having degree at least $3$.
\end{lm}

\subsection{The genus 3 case: the locus $\Sym^3(\calA_1)$}
\label{sec:genus3}
In genus 3 it is in fact known \cite{cmsurvey} that $I^{(3)}=\Sym^3(\calA_1)$ as a set (in the literature on the subject it is customary to write $\calA_1\times\calA_1\times\calA_1$ for the locus of such reducible ppav, while in fact it is really the symmetric product, and we try to make this distinction). The multiplicity is in fact equal to one as discussed in remark \ref{remred}.

We know from \cite[Th.~1.3]{grhu2} that  $\overline{G^{(3)}}= \overline{I^{(3)}}$.
We shall now compute the class of $\overline{G^{(3)}}$ and show that it coincides indeed with the class
$[\overline{\Sym^3(\calA_1)}]\in CH^*(\calA_3 ^{\op{Perf}})$ as computed in \cite{vdgeerchowa3}, and corrected in \cite{corr},
thus demonstrating the way our machinery works.

Indeed, if we use Lemma \ref{classlm}, the only new ingredient needed is an
expression for the terms with at least three $\delta$'s. We perform this computation in arbitrary genus for further use. To this end we first compute
$$
 \left(\sum\limits_{n\in Z_m}\de_n\right)^3=\sum\limits_{a\in Z_m}\de_a^3+3\sum\limits_{a\ne b\in Z_m}\de_a^2\de_b+6\sum\limits_{a<b<c \neq a \in Z_m}\de_a\de_b\de_c.
$$
By Proposition \ref{intersect} we have $2^{2g-2}$ sections $\overline{F_m}$ generically vanishing on any given boundary component $D_a$, $2^{2g-3}$ vanishing on any non-empty intersection $D_a\cap D_b$, and $2^{2g-4}$ vanishing on any non-empty intersection $D_a\cap D_b\cap D_c$ that is {\it global}, i.e.~over $\beta_3\subset\Sat$ (by part 1 of Proposition \ref{intersect} there do not exist any $\overline{F_m}$ generically vanishing on a {\it local} triple intersection of boundary divisors, over $\beta_2^0\subset\Sat$, which would correspond to the case $a+b+c=0$). By definition the locus $\beta_3\in \Perf(2)$ is the union of all global triple intersections, while $\sigma_3$ is the class of the union of {\it all} triple intersections, whether global or local. Altogether this yields
$$
 \sum\limits_{m{\rm\ odd}} {\overline p}_*\left(\sum\limits_{n\in Z_m}\delta_n^3\right)=
 2^{2g-2}{\overline p}_*(\sum\limits_a\de_b^3)+3\cdot2^{2g-3}{\overline p}_*\left(\sum\limits_{a\ne b} \de_a^2\de_b\right)
 +2^{2g-4}\cdot 6\frac{\beta_3}{2^3}
$$
where on the right-hand-side the summation is over all boundary
components labeled by $a,b,c\in (\ZZ/2)^{2g}\setminus 0$. The factor 6 in front of the last summand comes from the fact that every set $\{a,b,c\}$ appears 6 times in the summation, whereas the factor $8^3$ comes from the branching of order 8
along each boundary component.

We now need to deal with the combinatorics of the intersections of the boundary components. To this end, observe the identity
$$
{\overline p}_*\left((\sum\limits_a \de_a)(\sum\limits_{a< b} \de_a \de_b)\right)= \frac{\sigma_1\sigma_2}{2^3}={\overline p}_*\left(\sum\limits_{a\ne b}\de_a^2\de_b+ 3\sum\limits_{a<b<c}\de_a\de_b\de_c\right)
$$
We can thus deduce
$$
 {\overline p}_*\left(\sum\limits_{a\ne b}\de_a^2\de_b\right)=\frac{\sigma_1\sigma_2-3\sigma_3}{8},
$$
and similarly using
$$
 \frac{\sigma_1^3}{8}={\overline p}_*\left((\sum\limits_a \de_a)^3\right)={\overline p}_*\left(\sum\limits_a \de_a^3+3\sum\limits_{a\ne b}\de_a^2\de_b+6\sum\limits_{a<b<c}\de_a\de_b\de_c\right)
$$
we deduce
$$
 {\overline p}_*\left(\sum\limits_a \de_a^3\right)=\frac{\sigma_1^3-3(\sigma_1\sigma_2- 3\sigma_3)-6\sigma_3}{8}=\frac{\sigma_1^3-3\sigma_1\sigma_2+3\sigma_3}{8}.
$$
This finally yields, after grouping similar terms
\begin{lm}\label{3rdterm}
The term of the polynomial $P$ in Lemma \ref{classlm} involving products of three boundary divisors  equals
$$
 -\sum\limits_{m{\rm\ odd}} {\overline p}_*\left(\frac14\sum\limits_{n\in Z_m}\de_n\right)^3
 \sum\limits_{i=3}^g\binom{i}{3}\left(\frac{\lambda_1}{2}\right)^{i-3}\lambda_{g-i}
$$
$$
 =-2^{2g-12}\Big(2\sigma_1^3-3\sigma_1\sigma_2-3\sigma_3+3\beta_3\Big) \sum\limits_{i=3}^g\binom{i}{3}\left(\frac{\lambda_1}{2}\right)^{i-3}\lambda_{g-i}.
$$
\end{lm}
Combining this expression with the result of Lemma \ref{classlm} we finally get for
$g=3$
$$
  [\overline{G^{(3)}}]=28\left(\lambda_3+\lambda_2\frac{\lambda_1}{2} +\lambda_1\frac{\lambda_1^2}{4}+\frac{\lambda_1^3}{8}\right) -2\sigma_1\left(\lambda_2+2\lambda_1\frac{\lambda_1}{2}+3\frac{\lambda_1^2}{4}\right)
$$
$$  +\frac{1}{4}(\sigma_1^2-\sigma_2)\left(\lambda_1+3\frac{\lambda_1}{2}\right)
  -\frac{1}{32}(\sigma_1^3-3\sigma_1\sigma_2+3\sigma_3)-
  \frac{3}{64}(\sigma_1\sigma_2-3\sigma_3)-\frac{3}{64}\beta_3.
$$

The Chow ring $CH^*(\calA_3^{{\op {Perf}}})$ was computed in \cite{vdgeerchowa3}, and corrected in \cite{corr}, and in particular it was shown (as corrected in the erratum) that
$$
  [\overline{\Sym^3(\calA_1)}]=-35\lambda_3+\frac{35}{2} \lambda_1^3-\frac{25}{4}\lambda_1^2\sigma_1+\frac{5}{8}\lambda_1\sigma_2 +\frac{5}{8}\lambda_1\sigma_1^2-\frac{1}{12}\sigma_1\sigma_2.
$$
Simplifying the expression for $[\overline{G^{(3)}}]$ above and using the relations determined in \cite{vdgeerchowa3}, shows at the end of a lengthy, but straightforward, computation that we in fact have
$$
 [\overline{G^{(3)}}]=[\overline{\Sym^3(\calA_1)}]
$$
as expected.

\section{Genus 4 case: the class of the locus $\calA_1\times\theta_{\rm null}^{(3)}$ is tautological}
We will now treat the first new case, that of $g=4$. As remarked above, for $g=4$ the locus $\overline{I^{(4)}}$ is simply the closure of the locus of products $\calA_1\times\theta_{\rm null}^{(3)}$. Unlike the genus 3 case, $CH^*(\calA_4)$ and
$CH^*({\mathcal A}^{\operatorname {Perf}}_4)$
are not known (this is currently under investigation, see \cite{huto}), and in particular the classes of the locus $\calA_1\times\theta_{\rm null}^{(3)}$ and its closure are not known. However, the projections to the tautological ring of these classes can be computed, providing a good consistency check for theorem \ref{thmtaut}, which for $g=4$ gives the formula
$$
 [I^{(4)}]=45\lambda_1^4.
$$

In this section we first compute $[\calA_1\times\theta_{\rm null}^{(3)}]^{\op{taut}}$ directly geometrically, thus confirming our results, and then proceed to use theorem \ref{class} to obtain a complete expression for this class --- and not only for its tautological part.

\smallskip
We start by considering the locus of products $\calA_1\times\calA_3$.
\begin{lm}\label{gen4prod}
The projection of the class $[\overline{\calA_1\times\calA_3}]\in CH^*_\QQ(\calA_4^{\op{Perf}})$ to the tautological ring
is given by
$$
 [\overline{\calA_1\times\calA_3}]^{\op{taut}}= 20\lambda_3
$$
\end{lm}
\begin{proof}
We note that by \cite{vdgeercycles} there are only two degree 3 classes in the tautological ring, $\lambda_1^3$ and $\lambda_3$, and thus we know that we must have an expression of the form
$$
 [\overline{\calA_1\times\calA_3}]^{\op{taut}}= A\lambda_1^3+B\lambda_3,
$$
and it remains to determine the coefficients $A$ and $B$. This can be done by pairing with the complementary
dimension tautological classes. Indeed, by using the Hirzebruch-Mumford proportionality
principle \cite{mumhirz}, \cite{vdgeercycles}, these can be computed (we got the number for the top self-intersection number of $\lambda_1$ from \cite{vdgeercycles}, and used the relation $(1+\lambda_1+\ldots+\lambda_g) (1-\lambda_1+\ldots (-1)^g\lambda_g)=1$ from there to obtain $\lambda_3^2=\lambda_1^3\lambda_3-\lambda_1^6/8$ and $\lambda_1^5\lambda_3=7\lambda_1^8/48$)
$$
 \langle\lambda_1^3\cdot\lambda_1^7\rangle_{\calA_4^{\op{Perf}}}=\frac{1}{1814400};\quad \langle\lambda_3\cdot\lambda_1^7\rangle_{\calA_4^{\op{Perf}}}=\frac{7}{48}\cdot\frac{1}{1814400}
$$
$$
 \langle\lambda_1^3\cdot\lambda_1^4\lambda_3\rangle_{\calA_4^{\op{Perf}}}=\frac{7}{48}\cdot\frac{1}{1814400};\quad \langle\lambda_3\cdot\lambda_1^4\lambda_3\rangle_{\calA_4^{\op{Perf}}}=\frac{1}{48}\cdot\frac{1}{1814400}.
$$
We now need to compute the restrictions of the tautological classes to $\overline{\calA_1\times\calA_3}$ and compute the corresponding intersection numbers there.
Note that the perfect cone compactification is multiplicative \cite{shepherdbarron}, i.e.
$\overline{\calA_1\times\calA_3}= \calA_1^{\op{Perf}} \times \calA_3^{\op{Perf}}$.
Moreover, the Hodge bundle on the locus of products is the sum of the pullbacks of the Hodge bundles on the two factors,
so that we in particular have
\begin{equation} \label{restlambda1}
\lambda_1|_{\calA_1^{\op{Perf}} \times \calA_3^{\op{Perf}}}=1\times\lambda_1+\lambda_1\times 1
\end{equation}
and
\begin{equation} \label{restlambda3}
\lambda_3|_{\calA_1^{\op{Perf}} \times \calA_3^{\op{Perf}}}=1\times\lambda_3+\lambda_1\times \lambda_2.
\end{equation}
These expressions allow us to compute
$$\lambda_1^7[\overline{\calA_1\times\calA_3}]=7\langle\lambda_1\rangle_{\calA_1^{\op{Perf}}}
\cdot\langle\lambda_1^6\rangle_{\calA_3^{\op{Perf}}}=7\cdot\frac{1}{24}\cdot\frac{1}{181440}$$
and
$$\lambda_1^4\lambda_3[\overline{\calA_1\times\calA_3}]
=\langle \lambda_1\rangle_{\calA_1^{\op{Perf}}}
\cdot\langle \lambda_1^4\lambda_2\rangle_{\calA_3^{\op{Perf}}}+4\langle \lambda_1\rangle_{\calA_1^{\op{Perf}}}\cdot\langle \lambda_1^3\lambda_3\rangle_{\calA_3^{\op{Perf}}}
$$
$$
=\frac{1}{24}\cdot\frac12\cdot\frac{1}{181440}+4\frac{1}{24}\cdot\frac18\cdot\frac{1}{181440}=\frac{1}{24}
\cdot\frac{1}{181440}
$$
where for the last line we used $\lambda_2=\lambda_1^2/2$ and computed the $\lambda_1^3\lambda_3$ intersection number from $2\lambda_1\lambda_3=\lambda_2^2=\lambda_1^4/4$. Combining the above computations, we then get the following two equations for $A$ and $B$ (after factoring out the common factor of $1/1814400$):
$$
A+\frac{7}{48}B=10\cdot\frac{7}{24};\qquad \frac{7}{48}A+\frac{1}{48}B=10\cdot\frac{1}{24}
$$
and solving these gives the result of the proposition.
\end{proof}

\begin{cor}\label{cor:tn3}
The projection of the class $[\overline{\calA_1\times\theta_{\rm null}^{(3)}}]$ to the tautological ring (with rational coefficients)
is given by
$$
[\overline{\calA_1\times\theta_{\rm null}^{(3)}}]^{\op{taut}}=45\lambda_1^4.
$$
\end{cor}
\begin{rem}
To obtain the corollary we are basically going to intersect with a divisor, and deduce the tautological part of the product. It is tempting to say that this can be done in the tautological ring, but notice that a priori Chow is not a ring, and the projection from the Chow to the tautological ring cannot be a ring homomorphism.

Indeed, for example on $\overline{\calA_2}$ we have $\lambda_1^2\sigma_1=0$, and thus the projection of $\sigma_1$ to the tautological ring is equal to zero, while $\lambda_1\sigma_1^2\ne 0$, and thus the projection of $\sigma_1^2$ to the tautological ring is non-zero (in fact we have $\sigma_1^2=22\sigma_1\lambda_1-120\lambda_1^2$, see \cite{vdgeerchowa3}).
\end{rem}
\begin{proof}
{}From our previous proposition we know that
$$
[\overline{\calA_1\times \calA_3}]= 20\lambda_3 + [X]
$$
where the class $[X]$ has intersection $0$ with all monomials in the $\lambda_i$ of degree $7$.
The class of the theta-null divisor  $[\overline{\theta_{\rm null}^{(3)}}]$ in $CH^*(\calA_3^{\op{Perf}})$ equals
$18\lambda_1 - 2\sigma_1^{(3)}$ where $\sigma_1^{(3)}$ is the class of the boundary. We thus obtain
$$
18\lambda_1[\overline{\calA_1 \times \calA_3}] = [18\lambda_1|_{\calA_1^{\op{Perf}}} \times \calA_3^{\op{Perf}}]
+ [\overline{\calA_1\times\theta_{\rm null}^{(3)}}]
$$
$$
+ [\calA_1^{\op{Perf}} \times 2\sigma_1^{(3)}] +
18\lambda_1[X].
$$
We now deal with the summands term by term.
For the first term, we compute
$$
[18\lambda_1|_{\calA_1^{\op{Perf}}} \times \calA_3^{\op{Perf}}]=
\frac32[\calA_0\times\calA_3^{\op{Perf}}]
$$
(keeping in mind the stackiness, so that $\lambda_1$ has degree $1/24$ on $\calA_1$ and $\calA_0$ has degree $1/2$).
The tautological projection of the class $[\calA_0\times\calA_3^{\op{Perf}}]$ was computed in \cite[Prop.~4.3]{vdgeercycles},
resp. \cite[Theorem 3.4]{{ekvdgcycles}}, and we thus get
$$
 \frac32[\calA_0\times\calA_3^{\op{Perf}}]^{\op{taut}}=\frac{3\lambda_4}{2\zeta(-7)}=360\lambda_1\lambda_3-45\lambda_1^4.
$$

Let us now compute the tautological part of $[\calA_1^{\op{Perf}} \times 2\sigma_1^{(3)}]$ (note that lemmas \ref{projtotaut} and \ref{projzero} do not work here).
Indeed, this is a codimension 4 locus, and since we have $\lambda_1^6=8\lambda_1^3\lambda_3-8\lambda_3^2$, to compute the tautological part we need to know the intersections with $\lambda_3^2$ and with $\lambda_1^3\lambda_3$. Using formulae
(\ref{restlambda1}) and (\ref{restlambda3}) and recalling from
\cite[Lemma~3.11]{vdgeercycles} that in $CH^*(\calA_3^{\op{Perf}})$ we have $\lambda_3\sigma_1=0$, we thus get
$$
 \lambda_3[\calA_1^{\op{Perf}} \times 2\sigma_1^{(3)}]=[\lambda_1|_{\calA_1^{\op{Perf}}}
 \times \lambda_1^2\sigma_1^{(3)}],
$$
from which the intersection numbers are
$$
 \lambda_3^2[\calA_1^{\op{Perf}} \times 2\sigma_1^{(3)}]=0
$$
since $\lambda_1^2$ is zero in dimension one on $\calA_1^{\op{Perf}}$ and
$$
 \lambda_1^3\lambda_3[\calA_1^{\op{Perf}} \times 2\sigma_1^{(3)}]=0
$$
since $\lambda_1^5\sigma_1=0$ on $\calA_3^{\op{Perf}}$. Thus $[\calA_1^{\op{Perf}} \times 2\sigma_1^{(3)}]^{\op{taut}}=0$.

Since $[X]$ lies in the orthogonal complement of the tautological ring the same holds for $\lambda_1[X]$ and the assertion follows.
\end{proof}

\medskip
The above computation matches the result of our theorem \ref{thmtaut} for $g=4$. We now proceed to compute the class completely.
\begin{prop} \label{prop:classgenus4}
The class of the locus is
$$
 [\overline{I^{(4)}}]=[\overline{\calA_1\times\theta_{\rm null}^{(3)}}]=
180\lambda_1\lambda_3+\frac{45}{2}\lambda_1^4-8\sigma_1\lambda_3-14\sigma_1\lambda_1^3+\frac72
\lambda_1^2\sigma_1^2
$$
$$-\frac72\lambda_1^2\sigma_2-\frac38\lambda_1\sigma_1^3+\frac{9}{16}
\lambda_1\sigma_1\sigma_2+\frac{9}{16}\lambda_1\sigma_3-\frac{9}{16}\lambda_1\beta_3+\frac{3}{64}Y
+\frac{1}{64}\sigma_4
$$
$$-\frac{1}{16}\sigma_1\sigma_3+\frac{3}{64}\sigma_1\beta_3+\frac{1}{64}\sigma_2^2
-\frac{1}{32}\sigma_1^2\sigma_2+\frac{1}{64}\sigma_1^4
 \quad\in CH^*(\calA_4^{Perf})
$$
where the class $Y$ is defined as
$$
 Y:=\sum\limits_{n_1<n_2<n_3<n_4; n_1+n_2+n_3+n_4=0} \de_{n_1}\de_{n_2}\de_{n_3}\de_{n_4}.
$$
\end{prop}

\begin{rem}
One can understand the locus with class $Y$ geometrically, similarly to the way the discriminant locus $\Delta=\sigma_3-\beta_3$ can be described, in $\beta_2^0$ as the class of the stratum of semiabelic varieties whose normalization has two components, each of which is a $\PP^2$ bundle over some ppav $B\in\calA_{g-2}$. Indeed,
by inspecting the table of perfect cones for all possible small codimension strata in $\Perf$ we see that $Y$ is in fact equal to the locus of semiabelic varieties of torus rank 3 (since the dimension of the linear span of $n_i$ is three), for which the normalization is a union of two $\PP^3$ bundles and a $F(2,2)$ bundle. To prove this, we just observe that the corresponding cone is generated by $x_1^2,x_2^2,x_3^2, (x_1+x_2+x_3)^2$, giving the relation among the $n_i$ as in the definition of $Y$. We shall discuss this in more detail in remark \ref{rem:geom}.
\end{rem}
\begin{rem}
It is a very appealing question to try to define and describe a suitably ``extended'' tautological ring of a toroidal compactification of $\calA_g$ (in particular, of the perfect cone compactification). In addition to the tautological classes $\lambda_i$, presumably such an extended ring would include the classes $\lambda_i,\sigma_i,\beta_i$. In view of the above computation, it would also be natural to include $Y$ and the classes of all strata in the boundary corresponding to various types of semi-abelic varieties. Note that the proof below shows that $Y$ does not lie in the ring generated by $\sigma_i$ and $\beta_i$, but it could lie in the ring generated by these together with $\lambda_i$ --- we were unable to determine this.
\end{rem}

\begin{proof}[Proof of proposition \ref{prop:classgenus4}]
Note that one a priori technical difficulty is that $\calA_4^{\op{Perf}}$ is singular. However, it is only singular in codimension 10 (at the one point corresponding to the only non basic cone in the perfect cone compactification of genus $4$, namely the second perfect cone, see \cite{husaA4}).
The computations we do are for classes of codimension at most 4, and thus the singularity does not matter for these computations.

Indeed, applying the formulas from lemmas \ref{classlm} and \ref{3rdterm} and using the relations among the $\lambda$ classes from \cite{vdgeercycles} we get
$$
 [G^{(4)}]=180\lambda_1\lambda_3+\frac{45}{2}\lambda_1^4-8\sigma_1\lambda_1-14\sigma_1
 \lambda_1^3+\frac72\lambda_1^2\sigma_1^2-\frac72\lambda_1^2\sigma_2
$$
$$
 -\frac3{16}\lambda_1(2\sigma_1^3-3\sigma_1\sigma_2-3\sigma_3+3\beta_3)
$$
$$
 +\frac{1}{4^4}\sum\limits_{m{\rm\ odd}}{\bar p}_*\Big(24\sum\limits_{n_1<n_2<n_3<n_4\in Z_m}\delta_{n_1}\delta_{n_2}\delta_{n_3}\delta_{n_4}+12\sum\limits_{n_1;n_2<n_3\in Z_m}\delta_{n_1}^2\delta_{n_2}\delta_{n_3}
$$
$$
 +6\sum\limits_{n_1<n_2\in Z_m}\delta_{n_1}^2\delta_{n_2}^2
 +4\sum\limits_{n_1; n_2\in Z_m}\delta_{n_1}^3\delta_{n_2}
 +\sum\limits_{n\in Z_m}\delta_n^4\Big).
$$
where all the sums are taken over all $n_i\in Z_m$ that are pairwise distinct, and the coefficient of some monomial $\delta_{n_1}^{a_1}\cdot\ldots\cdot\de_{n_k}^{a_k}$ is equal to the binomial coefficient $\binom{\sum a_i}{a_1, \ldots, a_k}$

We now need to go through the possible combinatorics of the intersections. Note that any intersection $D_a\cap D_b$ for $a\ne b$, if non-empty, ``lives over $\beta_2$'' (which is to say that it intersects $\beta_2^0$ and does not intersect $\beta_1^0$), while for triple intersections we know that $D_a\cap D_b\cap D_{a+b}$ lives over $\beta_2$, and all other triple intersections live over $\beta_3$. Similarly the quadruple intersections that live over $\beta_3$ are $D_a\cap D_b\cap D_c\cap D_{a+b}$ and $D_a\cap D_b\cap D_c\cap D_{a+b+c}$, while all the other ones live over $\beta_4$, and in fact the other intersections together give $\beta_4$.

Using Proposition \ref{intersect}, and in particular that we cannot have $n_1+n_2+n_3=0$ for elements of $Z_m$, we thus compute for arbitrary $g$ (notice that the sums on the right are now over all $n_i$, not just those in $Z_m$; we suppress the fact that $n_i$ are ordered, and use $n$ for the index that is not ordered relative to $n_i$):
$$
\sum_m\sum\limits_{n_1<n_2<n_3<n_4\in Z_m}\!\!\!\!\delta_{n_1}\delta_{n_2}\delta_{n_3}\delta_{n_4}
=2^{2g-4}\!\!\!\sum\limits_{n_1 +n_2+n_3+n_4=0}\!\!\!\delta_{n_1}\delta_{n_2}\delta_{n_3}\delta_{n_4}
+2^{2g-5}\beta_4;
$$

$$
\sum_m\sum\limits_{n_1;n_2<n_3\in Z_m} \delta_{n_1}^2\delta_{n_2}\delta_{n_3}
=
2^{2g-4}\sum_{n_1+n_2+n_3\ne0; n_2<n_3}\de^2_{n_1}\de_{n_2}\de_{n_3};
$$
$$
\sum_m\sum\limits_{n_1<n_2\in Z_m}\de_{n_1}^2\de_{n_2}^2=
2^{2g-3}\sum\de_{n_1}^2\de_{n_2}^2,
$$
$$
\sum_m\sum\limits_{n_1,n_2\in Z_m}\de_{n_1}^3\de_{n_2}=
2^{2g-3}\sum\de_{n_1}^3\de_{n_2}.
$$
We now try to express the right-hand-sides as linear combinations of the degree four classes in the algebra generated by $\sigma_1=\beta_1,\sigma_2=\beta_2,\sigma_3,\beta_3,\sigma_4$, the corresponding expressions for which are
$$
  \sigma_4-\beta_4=\sum\limits_{n_1 +n_2+n_3+n_4=0}\delta_{n_1}\delta_{n_2}\delta_{n_3}\delta_{n_4}+
  \sum\limits_{n_1+n_2+n_3=0;n}\de_{n_1}\de_{n_2}\de_{n_3}\de_{n};
$$
$$
 \sigma_1\beta_3=4\beta_4+ 4\sum\limits_{n_1 +n_2+n_3+n_4=0}\delta_{n_1}\delta_{n_2}\delta_{n_3}\delta_{n_4}
 +3\sum\limits_{n_1+n_2+n_3=0;n}\de_{n_1}\de_{n_2}\de_{n_3}\de_{n}
$$
$$
 \hskip7cm+\sum\limits_{n_1+n_2+n_3\ne0; n_2<n_3}\de^2_{n_1}\de_{n_2}\de_{n_3},
$$
$$
 \sigma_1(\sigma_3-\beta_3)=\sum\limits_{n_1+n_2+n_3=0;n}\de_{n_1}\de_{n_2}\de_{n_3}\de_{n}
 +\sum\de^2_{n_1+n_2}\de_{n_1}\de_{n_2}
$$
for the terms involving some combinatorics. Notice that on the right-hand-side we have 4
different unknown terms, and thus we would not be able to express all of them in terms of these classes, so that the resulting expression will have to involve $Y$.

For the rest of the terms there is no combinatorics of the intersections involved, and we compute
$$
 \sigma_2^2=6\sum\limits_{\rm all}\de_{n_1}\de_{n_2}\de_{n_3}\de_{n_4}+2\sum\limits_{\rm all}\de_{n_1}^2\de_{n_2}\de_{n_3}+\sum\limits_{\rm all}\de_{n_1}^2\de_{n_2}^2;
$$
$$
 \sigma_1^2\sigma_2=12\sum\limits_{\rm all}\de_{n_1}\de_{n_2}\de_{n_3}\de_{n_4}+5\sum\limits_{\rm all}\de_{n_1}^2\de_{n_2}\de_{n_3}+2\sum\limits_{\rm all}\de_{n_1}^2\de_{n_2}^2
+\sum\limits_{\rm all}\de_{n_1}^3\de_{n_2};
$$
$$
\sigma_1^4=24\sum\limits_{\rm all}\de_{n_1}\de_{n_2}\de_{n_3}\de_{n_4}+12\sum\limits_{\rm all}\de_{n_1}^2\de_{n_2}\de_{n_3}
+6\sum\limits_{\rm all}\de_{n_1}^2\de_{n_2}^2;
$$
$$
+4\sum\limits_{\rm all}\de_{n_1}^3\de_{n_2}+\sum\limits_{\rm all}\de_n^4
$$
for the remaining terms (in the last three expressions we sum over all possible dimensions of linear spans of $\lbrace n_i\rbrace$ on the right).

Combining all of these allows us to first express all symmetric polynomials in $D_n$ in terms
of the standard symmetric polynomials,
and then we express all the terms involving various combinatorics of the intersections in terms of these classes and $Y$, obtaining the expression as claimed.
\end{proof}
For further use we record the result of the above computation for arbitrary $g$.
\begin{lm}
For arbitrary $g$ the order 4 term in proposition \ref{classlm} is equal to
$$
\frac{2^{2g-2}}{8^4}\left(6Y+3\beta_4+3(\sigma_1\beta_3-\beta_4-Y-3\sigma_4)
+3(\sigma_2^2-2\sigma_1\sigma_3+2\sigma_4)\right.
$$
$$\left.+2(\sigma_1^2\sigma_2-2\sigma_2^2-\sigma_1\sigma_3+4\sigma_4)
+\sigma_1^4-4\sigma_1^2\sigma_2+2\sigma_2^2+4\sigma_1\sigma_3-4\sigma_4
\right)
$$
$$%
=2^{2g-14}\left(\sigma_4-4\sigma_3\sigma_1+3Y+3\sigma_1\beta_3+\sigma_2^2-2\sigma_2\sigma_1^2 +\sigma_1^4\right)
$$
times the corresponding polynomial in $\lambda$ classes.
\end{lm}

\section{Genus 5: the locus of intermediate Jacobians}
In this section we finally compute the class in the Chow ring of the locus $I^{(5)}$, that is of the locus of intermediate Jacobians of cubic threefolds together with the locus of products $\calA_1\times\theta_{\rm null}^{(4)}$. We then compute the projection to the tautological ring of each of these two loci.
The computation is similar to the one in the previous section, but much more involved, with more new classes appearing.
Again we shall first work on the smooth part of $\calA_5^{\op{perf}}$. In this case we observe that the codimension of the singular locus is at least $6$, which is sufficient for our purpose:
away from $\beta_5$ we know it to be $10$, whereas on $\beta_5$ there is only one stratum of
codimension $5$ and this corresponds to the principal cone, which is basic.

\begin{prop}\label{g5class}
We have in $CH^*(\calA_5^{\op{Perf}})$ the formula
$$
 [\overline {IJ}]+[\overline{\calA_1\times\theta_{\rm null}^{(4)}}]=
 496\lambda_5+372\lambda_3\lambda_1^2+\frac{93}{2}\lambda_1^5-64\lambda_3\lambda_1\sigma_1-
 34\lambda_1^4\sigma_1
$$
$$
 +(4\lambda_3+14\lambda_1^3)(\sigma_1^2-\sigma_2)
 +\frac{5}{4}\lambda_1^2\left(3\sigma_3+3\sigma_2\sigma_1-2\sigma_1^3\right)
$$
$$
+\frac{7}{32}\lambda_1(\sigma_4-4\sigma_3\sigma_1+3Y+3\beta_3\sigma_1+\sigma_2^2-2\sigma_2
 \sigma_1^2+\sigma_1^4)
$$
$$
 -\frac{1}{256}\left(-95\sigma_5-30\beta_5-45A_2-30A_3-15A_4+15C_1+10D_1\right.
$$
$$
 \left.+45\sigma_4\sigma_1+15\beta_4\sigma_1+30Y\sigma_1
 +5\sigma_3\sigma_2-15\sigma_3\sigma_1^2+5\sigma_2^2\sigma_1-5\sigma_2\sigma_1^3+2\sigma_1^5\right)
$$
where the classes $A_2,A_3,A_4,C_1,D_1$ are defined below, and we refer to remark \ref{rem:geom} for a geometric interpretation of the classes $A_2,A_3,A_4$.
\end{prop}
\begin{proof}
We shall again use formula (\ref{Gclass2}).
Similar to the above computations, the new term we need to deal with here is
$\sum\limits_m (\sum\limits_{n\in Z_m}\de_n)^5$. Thus we will need to handle all products of five $\de_n$, with indices distinct or equal, and with all possible combinatorics of the intersections. The computations will be still more involved, and we proceed systematically, first listing all the terms that could appear. Recall that by convention the indices $n_1,\ldots,n_k$, and separately the indices $m_1,\ldots,m_\ell$ are ordered, and that all the linear relations among indices are stated under the sums. We then have the following types of products of 5 boundary components:
$$
A_1:=\sum \de_{n_1}\de_{n_2}\de_{n_3}\de_{n_4}\de_{n_5};\quad
A_2:=\sum\limits_{n_1+n_2+n_3+n_4+n_5=0} \de_{n_1}\de_{n_2}\de_{n_3}\de_{n_4}\de_{n_5};
$$
$$
A_3:=\sum\limits_{n_1+n_2+n_3+n_4=0} \de_{n_1}\de_{n_2}\de_{n_3}\de_{n_4}\de_{n};\quad
A_4:=\sum\limits_{n_1+n_2+n_3=0} \de_{n_1}\de_{n_2}\de_{n_3}\de_{m_1}\de_{m_2};
$$
$$
A_5:=\sum\limits_{n+n_1+n_2=n+m_1+m_2=0} \de_{n}\de_{n_1}\de_{n_2}\de_{m_1}\de_{m_2};
$$
$$
B_1:=\sum\de^2_{n}\de_{n_1}\de_{n_2}\de_{n_3};\quad
B_2:=\sum\limits_{n+n_1+n_2+n_3=0}\de^2_{n}\de_{n_1}\de_{n_2}\de_{n_3};
$$
$$
B_3:=\sum\limits_{n+n_1+n_2=0}\de^2_{n}\de_{n_1}\de_{n_2}\de_m;\quad
B_4:=\sum\limits_{n_1+n_2+n_3=0}\de^2_{n}\de_{n_1}\de_{n_2}\de_{n_3};
$$
$$
C_1:=\sum \de^2_{n_1}\de^2_{n_2}\de_{n};\quad
C_2:=\sum\limits_{n_1+n_2+n=0} \de^2_{n_1}\de^2_{n_2}\de_{n};
$$
$$
D_1:=\sum\de^3_n\de_{n_1}\de_{n_2};\quad
D_2:=\sum\limits_{n+n_1+n_2=0} \de^3_n\de_{n_1}\de_{n_2};
$$
$$
E:=\sum \de_n^3\de_k^2;\quad F:=\sum\de_n^4\de_k;\quad G:=\sum\de_n^5
$$
\begin{rem}\label{rem:geom}
Some of these loci --- the ones where each $\delta_n$ appears to power one --- have geometric interpretations.
Indeed, recall that the strata of a toroidal
compactification correspond to orbits of cones in the corresponding fan. The perfect cone compactification $\Perf$ is given,
as the name indicates, by the perfect cone decomposition of $\Sym^{2}_{\geq 0}(\RR^g)$ and $\Vor$ is given by the second Voronoi decomposition. For genus $g\leq 3$ these two toroidal compactifications coincide, whereas for $g=4,5$, but
not in general, the second Voronoi decomposition is a refinement of the perfect cone decomposition. In other words,
$\Vor$ is a blow-up of $\Perf$ for $g=4,5$ and by inspection of the decompositions one can see that the center of the
blow-up $\Vor \to \Perf$ has codimension $ > 5$. We also recall that, due to the moduli interpretation of $\Vor$ (see \cite{alexeev}), the strata
of this toroidal compactification have an interpretation in terms of polarized semi-abelic varieties.
In \cite[section 3]{grhu2} we have enumerated the relevant
cones of small codimension and described the corresponding semi-abelic varieties. We can use this description
to give a geometric interpretation for all of the above loci where each $\delta_n$ appears to power one.

First, we note that $A_1=\beta_5$. Then observe that $A_2$ is the class of the locus of semi-abelic varieties of torus rank 4 the normalization of which is a union of two $\PP^4$ bundles and one $X$ bundle (the toric polytope corresponding to $X$ is the 4-dimensional cube with two simplices removed) --- the corresponding cone is generated by $x_1^2,x_2^2,x_3^2,x_4^2,(x_1+x_2+x_3+x_4)^2$ and has the correct linear dependency. The class $A_3$ is the class of the stratum where the corresponding cone is generated by
$x_1^2,x_2^2,x_3^2,x_4^2,(x_1+x_2+x_3)^2$, and the normalization of the semiabelic variety is the union of two $\PP^1\times\PP^3$ bundles and one $\PP^1\times F(2,2)$ bundle. Similarly, $A_4$ corresponds to the stratum where the corresponding cone is $x_1^2,x_2^2,x_3^2,x_4^2,(x_1+x_2)^2$. The normalization of such a semiabelic variety is two copies of a $\PP^1\times\PP^1\times\PP^2$ bundle. Finally, $A_5$ corresponds to the cone $x_1^2,x_2^2,x_3^2,(x_1+x_2)^2,(x_1+x_3)^2$, and the normalization of the corresponding semiabelic varieties consists of two copies of a $\PP^3$ bundle and two copies of an $F_2$ (singular cone over $\PP^1\times\PP^1$) bundle.
\end{rem}

\smallskip
Denoting by $A$ the sum of all $A_i$, etc., we now express all degree 5 elements of the ring generated by $\sigma_i$ and $\beta_i$ in terms of the classes above. The expressions we get for polynomials in $\sigma_i$ are standard expressions in the algebra of symmetric polynomials, while in expressions involving $\beta_i$ the combinatorics of the indices plays a role --- some summands are missing, and the coefficients differ depending on the combinatorics. Indeed we have
$$
 \sigma_5=A;\quad \sigma_1\sigma_4=5A+B;\quad\sigma_2\sigma_3=10A+3B+C;
$$
$$
 \sigma_1^2\sigma_3=20A+7B+2C+D;\quad \sigma_1\sigma_2^2=30A+12B+5C+2D+E;
$$
$$
 \sigma_1^3\sigma_2=60A+27B+12C+7D+3E+F;
$$
$$
\sigma_1^5=120A+60B+30C+20D+10E+5F+G
$$
for the symmetric polynomials, and also
$$
 \beta_5=A_1;\quad \sigma_1\beta_4=5A_1+5A_2+4A_3+3A_4+B_1;
$$
$$
 \sigma_2\beta_3=10A_1+10A_2+10A_3+9A_4+8A_5+3B_1+3B_2+2B_3+3B_4+C_1;
$$
$$
 \sigma_1^2\beta_3=20A_1+20A_2+20A_3+18A_4+16A_5+7B_1+7B_2+5B_3+6B_4+2C_1+D_1.
$$
We want to use these identities to express all these unknown classes in terms of the smallest possible number of unknowns. We first note that from symmetric polynomial expressions we get expressions for $A,B,C,D,E,F,G$ in the ring generated by $\sigma_i$ (and thus writing down $A=A_1+A_2+A_3+A_4+A_5$ gives a non-trivial identity among the codimension 5 boundary strata --- which we then use to eliminate $A_5$ from the formulas). {}From the expression for $\sigma_1\beta_4$ we can express $B_1$ in terms of the geometrically defined classes $A_i$. We can further note that $Y\sigma_1=A_3+A_5+B_2$, which provides a geometric description of $B_2$.

\smallskip
Using proposition \ref{intersect}, we can now deal with the fifth order term appearing in the class computation in Theorem \ref{class}. Indeed, we note that no terms where there is a linear relation of odd length occur, and thus compute
$$
\frac{1}{4^5} p_*\sum\limits_{m{\rm\ odd}}\left(\sum\limits_{n\in Z_m}\de_n\right)^5=
 \frac{1}{8^5}\left(120(2^{2g-6}A_1+2^{2g-5}A_3)\right.
$$
$$
\left.
+60(2^{2g-5}B_1+2^{2g-4}B_2)
+30(2^{2g-4}C_1)+20(2^{2g-4}D_1)\right.
$$
$$
\left.
+10\cdot2^{2g-3}E+5\cdot2^{2g-3}F+2^{2g-2}G\right)
$$
$$
 =2^{2g-18}\left(15A_1+30A_3+15B_1+30B_2+15C_1+10D_1 \right.
$$
$$
\left.
+10E+5F+2G\right).
$$
To obtain an expression for this class in terms of the classes previously introduced, we note that $B_1$ and $B_2$ can be expressed in terms of the $A_i$ and $Y\sigma_1$, that we have eliminated $A_5$, but that there is no way to avoid using $C_1$ and $D_1$. Doing this computation in Maple yields the result of the proposition. Of course using the formulas for $\sigma_2\beta_3$ and $\sigma_1^2\beta_3$, we can express $C_1$ and $D_1$ in terms of $B_3$, $B_4$, and the other classes, but it is not clear why this would be a better expression.
\end{proof}

The proposition above computes the class of the locus $[I^{(5)}]$. To be able to compute the class of the locus of intermediate Jacobians, we need to compute the class of the other component. This is accomplished similarly to the computations in the previous section.
\begin{prop}
The projection of the class $[\overline{\calA_1\times\calA_4}]\in CH_\QQ^*(\calA_5^{\op{Perf}})$ to the tautological ring
is given by
$$
 [\overline{\calA_1\times\calA_4}]^{\op{taut}}= -\frac{11}{8}\lambda_1^4+11\lambda_1\lambda_3.
$$
In particular, if the locus of products $[\calA_1\times\calA_4]\in CH^*(\calA_5)$ is tautological, then $[\calA_1\times\calA_4]= -\frac{11}{8}\lambda_1^4+11\lambda_1\lambda_3$.
\end{prop}
\begin{proof}
The proof proceeds similarly to the proof of lemma \ref{gen4prod} for genus 4. {}From van der Geer's defining relation \cite[(1)]{vdgeercycles} of the tautological ring we get the following explicit identities in $CH^*(\calA_5^{\op{Perf}})$:
$$
 \lambda_3^2=\lambda_1^3\lambda_3-\frac{\lambda_1^6}{8}-2\lambda_1\lambda_5;\qquad
 \lambda_1^5\lambda_3= \frac{7}{48}\lambda_1^8+\frac{8\lambda_5}{3}(\lambda_1^3+\lambda_3);\qquad \lambda_5^2=0.
$$
Multiplying the first of these by $\lambda_1^5$ and the second by $\lambda_3$ and equating the results further gives
$$
 \lambda_1^3\lambda_3\lambda_5=\frac{1}{5}\lambda_1^6\lambda_5-\frac{1}{7040}\lambda_1^{11};
$$
multiplying this one by $\lambda_1^2$ and the second by $\lambda_5$ and equating the results finally yields
$$
 \lambda_1^8\lambda_5=\frac{3}{1144}\lambda_1^{13}.
$$
Using these relations, we can compute explicitly all the top intersection numbers in the ring generated by the $\lambda_i$ classes (of course, these also follow from the Hirzebruch-Mumford proportionality, but we could not find an easily available reference for intersection numbers on the symplectic Grassmannian).

Indeed, from the Hirzebruch-Mumford proportionality (see \cite{vdgeercycles} for the formula) we get $$\langle\lambda_1^{15}\rangle_{\calA_5^{\op{Perf}}}=\frac{13}{16329600}.$$
Using the above relations we compute step by step
$$
 \langle\lambda_1^{12}\lambda_3\rangle_{\calA_5^{\op{Perf}}}= \langle\frac83\lambda_1^7\lambda_3\lambda_5+\frac83\lambda_1^{10}\lambda_5+\frac{7}{48}\lambda_1^{15}
 \rangle_{\calA_5^{\op{Perf}}}
$$
$$
 =\langle\frac{16}{5}\lambda_1^{10}\lambda_5+\frac{8}{55}\lambda_1^{15}\rangle_{\calA_5^{\op{Perf}}}
$$
$$
=\frac{2}{13}\langle\lambda_1^{15}\rangle_{\calA_5^{\op{Perf}}}=\frac{2}{16329600}
$$
and then
$$
  \langle\lambda_1^{9}\lambda_3^2\rangle_{\calA_5^{\op{Perf}}}=
  \langle\lambda_1^{12}\lambda_3-2\lambda_1^{10}\lambda_5-\frac18\lambda_1^{15}\rangle_{\calA_5^{\op{Perf}}}
$$
$$
=
 \left(\frac{2}{13}-\frac{3}{572}-\frac18\right) \langle\lambda_1^{15}\rangle_{\calA_5^{\op{Perf}}}
  =\frac{1}{53222400}.
$$

{}From the fact that the Hodge bundle on $\calA_5$ restricts to the sum of Hodge bundles on the factors of a decomposable ppav, we get, as in lemma \ref{gen4prod},
$$\lambda_1^{11}[\overline{\calA_1\times\calA_4}]=11\langle\lambda_1\rangle_{\calA_1^{\op{Perf}}}
\cdot\langle\lambda_1^{10}\rangle_{\calA_4^{\op{Perf}}}=\frac{11}{24}\cdot\frac{1}{1814400}$$
and
$$\lambda_1^8\lambda_3[\overline{\calA_1\times\calA_4}]=\langle \lambda_1\rangle_{\calA_1^{\op{Perf}}}
\cdot\langle \lambda_1^8\lambda_2\rangle_{\calA_4^{\op{Perf}}}+8\langle \lambda_1\rangle_{\calA_1^{\op{Perf}}}\cdot\langle \lambda_1^7\lambda_3\rangle_{\calA_4^{\op{Perf}}}
$$
$$
=\frac{1}{24}\cdot\frac12\cdot\frac{1}{1814400}+8\cdot\frac{1}{24}\cdot\frac{7}{48}\cdot\frac{1}{1814400}
=\frac{5}{72}
\cdot\frac{1}{1814400}
$$
where we have used the computation of $\langle \lambda_1^7\lambda_3\rangle_{\calA_4^{\op{Perf}}}$ from lemma \ref{gen4prod}.

Thus if we have $[\overline{\calA_1\times\calA_4}]=A\lambda_1^4+B\lambda_1\lambda_3$ in the tautological ring, we can compute the coefficients from the above relations, and the result is as stated.
\end{proof}
Similarly to the $g=4$ case, we can then compute the projection of the class of the theta-null divisor to
the tautological ring.
\begin{cor}
The projection of the class $[\overline{\calA_1\times\theta_{\rm null}^{(4)}}]$ to the tautological ring with rational coefficients is given by
$$
[\overline{\calA_1\times\theta_{\rm null}^{(4)}}]^{\op{taut}}=187(-\lambda_1^5/2+4\lambda_1^2\lambda_3 - 4\lambda_5).
$$
\end{cor}
\begin{proof}
Indeed, similarly to the proof of corollary \ref{cor:tn3}, let $[X]\in CH^4(\calA_5)$ be the part of the class of $[\overline{\calA_1\times\calA_4}]$ orthogonal to the tautological ring.
Since the class of $[\overline{\theta_{\rm null}^{(4)}}]$ is equal to $68\lambda_1-8\sigma_1^{(4)}$, we compute
$$
 68\lambda_1[\overline{\calA_1\times\calA_4}]=
[68\lambda_1|_{\calA_1^{\op{Perf}}}\times\calA_3^{\op{Perf}}]+[\overline{\calA_1\times\theta_{\rm null}^{(4)}}]
+68\lambda_1[X]+8[\calA_1^{\op{Perf}}\times \sigma_1^{(4)}].
$$
The class $\lambda_1[X]$ is orthogonal to the tautological ring, and we need to argue that so is the last class. Indeed, it is codimension 5 in $\calA_5^{\op{Perf}}$, but by the restriction property of the Hodge bundle it is enough to argue that $\sigma_1^{(4)}$ is orthogonal to the tautological ring of $\calA_4^{\op{Perf}}$, which follows from lemma \ref{projtotaut}.

Thus using the above proposition
together with the formula $[\calA_0 \times \calA_4^{\op{Perf}}]^{\op{taut}} = -{\lambda_5}/{\zeta(-9)}$ from \cite[prop.~4.3]{vdgeercycles} (and recalling again that $\deg\calA_0=1/2$, as of a stack, and thus $\lambda_1|_{\calA_0}=1/12[\calA_0]$),
we obtain
$$
-2\cdot\frac{68\lambda_5}{24\zeta(-9)}+[\overline{\calA_1\times\theta_{\rm null}^{(4)}}]^{\op{taut}}
=187(-\lambda_1^5/2 +4\lambda_1^2\lambda_3).
$$
Recalling $\zeta(-9)=-1/132$ proves the claim.
\end{proof}
Combining this with theorem \ref{thmtaut} finally yields proposition \ref{IJtaut}, giving a formula for $[\overline{IJ}]^{\op{taut}}$

\section{Degeneration of intermediate Jacobians}
\label{degen}
As a conclusion we discuss the geometry of the boundary of the locus of intermediate Jacobians of cubic threefolds. While the moduli space of cubic threefolds has been constructed as a ball quotient by Allcock, Carlson, and Toledo \cite{alcato}, its boundary is not yet completely understood. {}From the point of view of intermediate Jacobians, Casalaina-Martin and Laza \cite[Theorem 1.1]{cmla} described the boundary of the closure $\overline{IJ^0}^{\op{Sat}}$ of the locus $IJ^0$ of intermediate Jacobians
in the Satake compactification $\calA_5^{\op {Sat}}$.
Indeed, they showed that this consists of three $9$-dimensional irreducible components, namely $\calJ_5^h$, the locus of Jacobians of hyperelliptic genus $5$ curves (contained in the locus of indecomposable ppav), $\calA_1\times(\calJ_4\cap\theta_{\rm null}^{(4)})$ (contained in the locus of decomposable ppav), and $\calJ_4$, contained in the boundary $\calA_4\subset\partial \calA_5^{\op {Sat}}$. We note also that $\calA_3\subset\partial \calA_5^{\op {Sat}}$ is contained in the boundary of $\calJ_4$, and thus is contained in the boundary of $\overline{IJ^0}$ in the Satake compactification.

The closure of $I^{(5)}$ in the partial toroidal compactification of $\calA_5$ was described by Salvati Manni and the first author in \cite{grsmconjectures}.

The boundary of the toroidal compactifications $\calA_5^{\op {Vor}}$ and $\calA_5^{\op {Perf}}$ is a $\QQ$-Cartier
divisor and hence the boundary of $\overline{IJ}$ and of $\overline{I^{(5)}}$ in these spaces is of
pure dimension $9$, in
particular every irreducible component has this dimension. It currently seems very hard to describe the boundary
in $\calA_5^{\op {Vor}}$ and hence we shall concentrate on the boundary of $\overline{IJ}$ and $\overline{I^{(5)}}$
in $\calA_5^{\op {Perf}}$. Note that it follows from the computation
\cite{grhu2} that every such boundary component intersects the partial compactification.

We shall now apply our results from \cite{grhu2} to enumerate the boundary components
of $\overline{IJ^0}$ and of $\overline{I^{(5)}}$ (in $\calA_5^{\op {Perf}}$) and describe
their generic points.
At this point we would like to recall the approach of Casalaina-Martin and Laza
for studying the boundary of $IJ^0$.
They use the fact that the intermediate Jacobian of a cubic threefold $X$ is the Prym variety
of a double cover of a plane quintic $C$. As $X$ acquires singularities, so does $C$, and it is possible to
have some control over the singularities of the plane quintic as well as the possible admissible covers, see
in particular \cite[section 5]{cmla}.
We note that our result implies in particular the results of \cite{cmla} for the components of the
boundary of $\overline{IJ^0}^{\op{Sat}}$.
For the other component $\calA_1\times \theta_{\rm null}^{(4)}$ of $I^{(5)}$, one sees easily that its boundary in the Satake compactification is equal to $\theta_{\rm null}^{(4)}$ (which contains $\calA_1\times\calA_3$).

A more precise result, describing the boundary of $I^{(5)}$ in the partial toroidal compactification, is \cite[Proposition 12]{grsmconjectures}, and this was also studied by us in much more detail in section 4 of \cite{grhu2}.
The result involves the global family of singularities of theta divisors $\calS_g\subset\calX_g$, and its projection to $\calA_g$, the locus of ppav whose theta divisor is singular, i.e.~the Andreotti-Mayer divisor $N_0$. We recall (see \cite{civdg1},\cite{grsmordertwo}) that the locus $\calS_g$ is equidimensional, of codimension $g+1$, and has three irreducible components: $\calS_{\rm dec}$ projects to $\calA_1\times\calA_{g-1}$, with $(g-2)$-dimensional fibers being the product of theta divisors; $\calS_{\rm null}$ projects to $\theta_{\rm null}$ generically one-to-one (and finitely), and the remaining component $\calS'$ projects to the component $N_0'$ of the Andreotti-Mayer divisor. In particular, for $g=4$ it projects onto $\calJ_4$.

Recall that the partial boundary $\beta_1^0$ is equal to the universal Kummer family $\calX_4/\pm 1$, a point on which we denote $(\tau,b)$ for $\tau\in\calA_4$. The result \cite[prop.~12]{grsmconjectures} is that
\begin{equation}\label{bdryI}
  \overline{I^{(5)}}\cap \beta_1^0=2_*(\calS)\cup\pi^{-1}(I^{(4)}),
\end{equation}
where $2_*$ denotes the multiplication by two on each ppav, and $\pi:\calX_4\to\calA_4$.. Thus the partial boundary $\overline{I^{(5)}}\cap \beta_1^0$ is the union of four loci:
$$
\begin{aligned}
 (I):=&\lbrace (\tau,b)\mid \tau\in\calJ_4, b\in2_*(\Sing\Theta_\tau)\rbrace,\\
 (II):=&\lbrace (\tau,b)\mid \tau\in\theta_{\rm null}^{(4)}, b=0\rbrace,\\
 (III):=&\lbrace (\tau,b)\mid \tau\in \calA_1\times \theta_{\rm null}^{(3)}, b\in2_*(\Theta_\tau)\rbrace,\\
 (IV):=&\lbrace (\tau,b)\mid \tau\in\calA_1\times\calA_3, b \in 2_*(\Sing\Theta_\tau)\rbrace,
\end{aligned}
$$
projecting respectively to $\calJ_4,\theta_{\rm null}^{(4)}$, $\calA_1\times \calJ_3^h$, and $\calA_1\times\calA_3$ in the boundary of the Satake compactification.

For $(I)$, we note that for a non-hyperelliptic Jacobian of a genus 4 curve the two singular points on its theta-null are the differences of the two $g^1_3$'s on the curve; the corresponding $b$ is thus (generically) unique, and the corresponding locus is irreducible.
For $(II)$, $b$ is really twice the corresponding two-torsion point, but this is zero, and since $\theta_{\rm null}^{(4)}$ is irreducible, the locus $(II)$ is irreducible.
For $(III)$ we used the fact that $I^{(4)}=\calA_1\times \theta_{\rm null}^{(3)}=\calA_1\times \calJ_3^h$; then we either have $b=2(m\times z)$ for $m$ being the (unique) odd two-torsion point on the elliptic curve and $z$ arbitrary (so $b\in \lbrace0\rbrace\times A_{\tau_3}$), or $b=2(x\times z)$, where $x$ is arbitrary, and $z$ lies on the theta divisor of the Jacobian of a hyperelliptic curve --- thus $(III)$ has two irreducible components.
Finally, in $(IV)$ the singular locus of the theta divisor is generically $m\times \Theta^{(3)}$, and thus $b=0\times z$, where $z$ lies on twice the theta divisor of the abelian threefold factor --- so the locus $(IV)$ is also irreducible.

The boundary of $\overline{\calA_1\times\theta_{\rm null}^{(4)}}$ in $\beta_1^0$ can be easily described geometrically: there are two possibilities depending on which factor degenerates. Indeed, if the elliptic curve degenerates to a rational nodal curve, we get the locus of trivial extensions over $\theta_{\rm null}^{(4)}$ i.e.~the component $(II)$. On the other hand, the boundary of the theta-null divisor in the partial compactification can be described scheme-theoretically similarly to the above (see \cite{grsmconjectures}): we have
\begin{equation}\label{bdrythetanull}
  \overline{\theta_{\rm null}^{(4)}}\cap \beta_1^0=2_*(\calT)\cup\pi^{-1}(\theta_{\rm null}^{(3)}),
\end{equation}
where $\calT\subset\calX_3$ denotes the universal theta divisor (recall that it is of course defined only up to translation by a two-torsion point), but is image under the multiplication by two is well-defined. Thus for this case for the partial boundary of $\overline{\calA_1\times\theta_{\rm null}^{(4)}}$ we get two more irreducible components: the locus where $\tau\in\calA_1\times\calA_3$, and $b$ is of the form $0\times z$ for $z$ on the theta divisor, i.e.~case $(IV)$, and the locus of $\tau\in\theta_{\rm null}^{(3)}$, with $b$ arbitrary --- in this case for $\overline{\calA_1\times\theta_{\rm null}^{(4)}}\cap\beta_1^0$ we get the first of the two irreducible components of case $(III)$.

{}From the results of \cite{cmla} it follows that $(II)$ does not lie in the boundary of $\overline{IJ}$ (since its projection to the Satake compactification does not); since $(I)$ and the second irreducible component of $(III)$ do not lie in the boundary of $\calA_1\times\theta_{\rm null}^{(4)}$, they must lie in the boundary of $\overline{IJ}$.
We claim that this describes completely the boundary of the locus of intermediate Jacobians in the partial compactification:
\begin{thm}
The boundary of $\overline{IJ}^{\op{Perf}}\subset{\calA_5^{\op {Perf}}}$ has exactly two irreducible
components: these are the closures of $(I)$ and of the second component of $(III)$ above, i.e.~the two irreducible components are the closures of the two components of
$$
 \begin{aligned}
  \overline{IJ}^{\op{Perf}}\cap\beta_1^0=&\lbrace (\tau,b)\mid \tau\in\calJ_4, b\in2_*(\Sing\Theta_\tau)\rbrace\\
  \cup &\lbrace (\tau,b)\mid \tau=\tau_1\times \tau_3, b\in A_{\tau_1}\times 2_*(\Theta_{\tau_3})\}
 \end{aligned}
$$
where $\tau_1\in\calA_1$, $\tau_3\in\theta_{\rm null}^{(3)}$, and $(A_{\tau_i},\Theta_{\tau_i})$ denote the corresponding ppav.
\end{thm}
\begin{proof}
{}From our explicit description of semiabelic theta divisors on all strata in $\Perf$ of codimension at most 5, from \cite{grhu2}, it follows, by dimension reason that the boundary of $\overline{IJ}^{\op{Perf}}\subset{\calA_5^{\op {Perf}}}$ cannot have any irreducible components contained in $\beta_2$, and thus to prove the theorem it remains to show that $(IV)$ and the first component of $(III)$ are not contained in $\overline{IJ}$. We remark that formula \ref{bdryI} describing the intersection of $\overline{I^{(5)}}$ with $\beta_1^0$ as the union of $2_*\calS\subset\calX_4$ (the universal locus of singularities of the theta divisor) and $\pi^{-1}(\overline{I^{(4)}})$
in fact holds  scheme-theoretically, as the proof in \cite{grsmconjectures} is by studying the defining equations.

It is known that $\calS\subset\calX_4$ is equidimensional, of dimension nine, and has three irreducible components:
$\calS=\calS_{\rm null}\cup\calS_{\rm dec}\cup\calS'$ (see \cite{debarredecomposes},\cite{grsmconjectures}), where $\calS_{\rm null}$ projects to $\theta_{\rm null}^{(4)}$, $\calS_{\rm dec}$ is equal to $(IV)$ above, and projects to $\calA_1\times\calA_3$, and $\calS'$, which projects to $\calJ_4$, and is equal to $(I)$ above. We recall
that all three of these components are reduced: for $\calS'$ and $\calS_{\rm null}$ this is immediate as they project generically finitely to their images in $\calA_4$, which are divisors there, while the fact that $\calS_{\rm dec}$ is also reduced seems to be well-known: directly this can be seen by computing the Jacobian matrix of the defining equations at a generic point, similarly to the proof of \cite[thm.~6]{grsmconjectures}. Thus the intersection $\overline{I^{(5)}}\cap\beta_1^0$ is reduced, and in particular the two irreducible components $\overline{\calA_1\times\theta_{\rm null}^{(4)}}$ and $\overline{IJ}$ of $\overline{I^{(5)}}$ must intersect $\beta_1^0$ in different component.
It follows in particular that the first irreducible component of $(III)$ and $(IV)$
lie only in $\overline{\calA_1\times\theta_{\rm null}^{(4)}}$, and are not contained in $\overline{IJ}$.
\end{proof}

\begin{rem}
We note that the second component of $(III)$ does not appear in the results of \cite{cmla} on the boundary of $\overline{IJ}^{\op{Sat}}$ in the Satake compactification, as its projection to the Satake compactification has $3$-dimensional fibers (while $(I)$ projects to the Satake compactification generically one-to-one), and thus its Satake image is higher codimension. It is of course a natural question to also study the closure $\overline{IJ}^{\op{Vor}}$ in the second Voronoi compactification of $\calA_5$, but at the moment this seems out of reach.
\end{rem}

Our more detailed results in \cite{grhu2} on the vanishing loci of $f_m$ allow one to describe explicitly the geometry of the boundary of $\overline{I^{(5)}}$ in other strata of the perfect cone compactification. However, these computations become considerably more involved, and we will investigate this subject in detail in a future work. Here we just give one such result to give a flavor of the geometry involved.

Indeed, the next natural stratum to consider beyond $\beta_1^0$ is torus rank two,
i.e.~the set $\beta_2^0$, which is the pre-image of $\calA_{g-2}$ under the projection of $\Perf$ to the Satake compactification.
Recall that there are then two possibilities: on the open part of $\beta_2^0$ the normalization of the semiabelic variety is a $\PP^1\times \PP^1$ bundle over a $(g-2)$-dimensional ppav (these are called standard degenerations in \cite{grhu2}), while over a codimension 3 stratum $\Delta\subset\beta_2^0$ ($\Delta$ is codimension 3 in $\Perf$) the normalization of the semiabelic variety consists of two $\PP^2$ bundles.

We give the description of $I^{(g)}\cap\Delta$ resulting from \cite{grhu2}; while for the standard rank two degenerations the description is easier, it is perhaps less significant geometrically, as it can be more easily understood as a sequence of two degenerations.

Recall that as a set $\Delta$ is the fiberwise square $\calX_{g-2}^{\times 2}$, points on which we denote $(\tau,b_1,b_2)$, for $\tau$ corresponding to a ppav $(A_\tau,\Theta_\tau)\in\calA_{g-2}$, with points $b_1,b_2\in A_\tau$. The geometry was first investigated by Ciliberto and van der Geer in \cite{civdg2}, and studied in detail in \cite[section 6]{grhu2}. The result is as follows
\begin{prop}
The intersection $\overline{I^{(g)}}\cap\Delta$ of the closure of the locus $I$ with the stratum $\Delta$ has the following (not necessarily irreducible) components
$$
\begin{aligned}
 (i):=&\lbrace (\tau,b_1,b_2)\mid x:=\frac{b_1-b_2}{2}\in\Sing\Theta_\tau, b_2+x\in \Theta_\tau \rbrace,\\
 (ii):=&\lbrace (\tau,b_1,b_2)\mid x:=\frac{-b_2}{2}\in\Sing\Theta_\tau, b_1+x\in\Theta_\tau\rbrace,\\
 (iii):=&\lbrace (\tau,b_1,b_2)\mid x:=\frac{-b_1}{2}\in\Sing\Theta_\tau, b_2+x\in\Theta_\tau\rbrace,\\
 (iv):=&\lbrace (\tau,b_1,b_2)\mid \tau\in I^{(g-2)} , b_1,b_2\in \Theta_\tau+m \rbrace.
 \end{aligned}
$$
In the first three cases, by dividing by two we mean that this $x$ can be chosen up to a two-torsion point, and in the last case $m$ is the odd two-torsion point where $\Theta_\tau$ has multiplicity three.
\end{prop}

The first three cases correspond to the case of singularities at one of the $\PP^1$ ``edges'' of the triangle polytope corresponding to $\PP^2$, and the fourth case corresponds to the singularity at the vertex of the triangle (all vertices are glued). For $g=5$, all of these cases can be described explicitly: in each of the cases $(i),(ii),(iii)$ we have two irreducible ($7$-dimensional) components of $\overline{I^{(5)}}\cap\Delta$, one projecting to $\theta_{\rm null}^{(3)}$ in the Satake, with $b_1-b_2$ (resp.~$b_2,b_1$) being zero, and the point $b_2$ (resp.~$b_1,b_2$) lying on the theta divisor, and the other irreducible component projecting onto $\calA_1\times\calA_2$, with the point $b_1-b_2$ (resp.~$b_2,b_1$) lying in $0\times C$ (where $C$ is the curve of which the abelian surface is the Jacobian), and the point $b_2$ (resp.~$b_1,b_2$)
lying on the theta divisor $(m\times A_{\tau_3})\cup (A_{\tau_1}\times C)$.
For case $(iv)$, we recall that $I^{(3)}=Sym^3(\calA_1)$, and $b_1$ and $b_2$ must then lie on the theta divisor of $E_1\times E_2\times E_3$, shifted by $m$, i.e.~on one of the three ``coordinate'' products $E_i\times E_j$ embedded in it (the dimension of this case is of course still $7=3+2+2$).

\begin{rem}
All except the first component of $(i)$ above are contained in the boundary of $\calA_1\times\theta_{\rm null}^{(4)}$: the first components of $(ii),(iii)$ are contained in the closure of $(II)$, the second components of $(i),(ii),(iii)$ --- in the closure of $(IV)$, and $(iv)$ --- in the closure of the first component of $(III)$.
Since in the first component of $(i)$ we have $b_1=b_2$, arbitrary lying on $\Theta_\tau$, and in particular neither $b_1$ nor $b_2$ has a zero component, this component does not lie in the boundary of $\calA_1\times\theta_{\rm null}^{(4)}$, and must lie in the boundary of $\overline{IJ}^{\op{Perf}}$. However, it is not clear to us whether $\overline{IJ}^{\op{Perf}}$ contains any of the remaining components of $\overline{I^{(5)}}\cap\Delta$ listed above. Note also that to list completely all components of the boundary of $\overline{I^(5)}$ in $\Perf$ of codimension up to 3, one would need to perform similar analysis for the standard rank two and three degenerations, and to ensure that there can be no further components of this codimension contained in $\beta_4$ --- our analysis in \cite{grhu2} does not quite guarantee that.

Thus a further detailed analysis of the boundary of the locus of intermediate Jacobians is warranted, and of great interest also for compactifying geometrically the moduli of cubic threefolds. We will further pursue this topic elsewhere.
\end{rem}

\end{document}